\documentclass[amssymb, 11pt]{amsart}
\usepackage{amsmath}
\usepackage{amsthm}
\usepackage{amsfonts}
\newcommand{\pp}{\mathbb{P}}
\newcommand{\ee}{\mathbb{E}}
\newcommand{\qb}[2]{{\left [{#1 \atop #2} \right
]}}
\newtheorem{thm}{Theorem}[section]
\newtheorem{theorem}[thm]{Theorem}

\newtheorem{cor}[thm]{Corollary}
\newtheorem{prop}[thm]{Proposition}

\newtheorem{lemma}[thm]{Lemma}

\newcommand{\num}{\refstepcounter{thm}}

\begin{document}

\title [Commutation relations and Markov chains] {Commutation relations
and Markov chains}

\author{Jason Fulman}
\address{Department of Mathematics\\
University of Southern California\\
Los Angeles, CA 90089}
\email{fulman@usc.edu}



\thanks{Submitted December 9, 2007; referee suggestions implemented on January 20, 2008.}

\begin{abstract} It is shown that the combinatorics of commutation relations is
well suited for analyzing the convergence rate of certain Markov chains.
Examples studied include random walk on irreducible representations, a
local random walk on partitions whose stationary distribution is the Ewens
distribution, and some birth-death chains.
\end{abstract}

\maketitle

\section{Introduction} \label{intro}

Stanley \cite{St2} introduced a class of partially ordered sets, which he
called differential posets, with many remarkable combinatorial and
algebraic properties. A basic tool in his theory was the use of two linear
transformations $U$ and $D$ on the vector space of linear combinations of
elements of $P$. If $x \in P$ then $Ux$ (respectively, $Dx$) is the sum of
all elements covering $x$ (respectively, which $x$ covers). For
differential posets one has the commutation relation $DU-UD=rI$ for some
positive integer $r$, and he exploited this to compute the spectrum and
eigenspaces (though typically not individual eigenvectors) of the operator
$UD$.

The primary purpose of this paper is to show that commutation relations
are useful not only for studying spectral properties, but also for
obtaining sharp Markov chain convergence rate results. We will need the
more general commutation relation (studied in Fomin's paper \cite{Fo})
\num \begin{equation} \label{gencom} D_{n+1}U_n = a_n U_{n-1} D_n + b_n
I_n,
\end{equation} for all $n$. In many of our examples the operators $U,D$
will not be Stanley's up and down operators but will be probabilistic in
nature and will involve certain weights.

There are several ways of quantifying the convergence rate of a Markov
chain $K$ to its stationary distribution $\pi$. These, together with other
probabilistic essentials, will be discussed in Section \ref{probback}. For
now we mention that the commutation relations (\ref{gencom}) will be
particularly useful for studying the maximal separation distance after $r$
steps, defined as \[ s^*(r) := \max_{x,y} \left[ 1 -
\frac{K^r(x,y)}{\pi(y)} \right], \] where $K^r(x,y)$ is the chance of
transitioning from $x$ to $y$ in $r$ steps. In general it can be quite a
subtle problem even to determine which $x,y$ attain the maximum in the
definition of $s^*(r)$. Our solution to this problem involves using the
commutation relations (\ref{gencom}) to write $\frac{K^r(x,y)}{\pi(y)}$ as
a sum of non-negative terms.

After determining which $x,y$ maximize $1 - \frac{K^r(x,y)}{\pi(y)}$,
there is still work to be done in analyzing the value of $s^*(r)$, and in
particular its asymptotic behavior. For several examples in this paper,
our method of writing $\frac{K^r(x,y)}{\pi(y)}$ as a sum of non-negative
terms will be well-suited for this. For all of the examples in this paper,
we do express $s^*(r)$ in terms of the distinct eigenvalues
$1,\lambda_1,\cdots,\lambda_d$ of $K$:
\begin{enumerate}
\item \[ s^*(r) = \sum_{i=1}^d \lambda_i^r \left[ \prod_{j \neq i}
\frac{1-\lambda_j}{\lambda_i-\lambda_j} \right].\]

\item \[ s^*(r) = \pp(T>r), \] where $T = \sum_{i=1}^d X_i$ and the $X_i$
are independent geometric random variables with probability of success
$1-\lambda_i$.
\end{enumerate} These relations are useful for studying convergence rates
and appeared earlier for certain one-dimensional problems (stochastically
monotone birth-death chains started at $0$) \cite{DF},\cite{DSa}, and in
\cite{F3} for a higher dimensional problem (random walk on irreducible
representations of $S_n$). The current paper provides further examples,
and revisits the results of \cite{F3} using commutation relations.

Section \ref{downup} reviews the concept of ``down-up'' Markov chains on
branching graphs and describes some main examples to be analyzed in this
paper. Aside from their intrinsic combinatorial interest, down-up chains
are very useful. They were crucially applied in \cite{F1}, \cite{F4} to
study asymptotics of characters of the symmetric group, and were recently
used in \cite{BO2},\cite{Pe} to construct interesting infinite dimensional
diffusions. (Actually \cite{BO2}, \cite{Pe} use ``up-down'' chains instead
of ``down-up'' chains; our methods apply to these too and it will be shown
that they have the same convergence rate asymptotics). Convergence rate
information about these chains is also potentially useful for proving
concentration inequalities for statistics of their stationary
distributions \cite{C}.

Section \ref{dfrep} adapts Stanley's work on differential posets to the
commutation relations (\ref{gencom}). These results are applied in Section
\ref{young} to study the down-up walk on the Young lattice. Here the
stationary distributions are the so called z-measures, studied in papers
of Kerov, Olshanski, Vershik, Borodin, and Okounkov (see \cite{KOV1},
\cite{KOV2}, \cite{BO2}, \cite{BO3}, \cite{O1} and the references
therein). In a certain limit these measures become the Plancherel measure
of the symmetric group, and we obtain new proofs of results in \cite{F3}.

Sections \ref{schur} analyzes down-up walk on the Schur lattice. We
explicitly diagonalize this random walk, and use this to study total
variation distance convergence rates. Similar ideas can be used to analyze
down-up walk on the Jack lattice (see the discussion at the end of Section
\ref{schur}). The arguments in Section \ref{schur} do not require the use
of commutation relations, though we do note some connections.

Section \ref{kingman} applies commutation relations to study down-up walk
on the Kingman lattice. Here the stationary distribution depends on two
parameters $\theta,\alpha$ and when $\alpha=0$ is the Ewens distribution
of population genetics. The down-up walk is more ``local" than the
traditionally studied random walks with this stationary distribution, such
as the random transposition walk when $\alpha=0,\theta=1$; this could be
useful for Stein's method. We show that the eigenvalues and separation
distance do not depend on the parameter $\alpha$, and prove order $n^2$
upper and lower bounds for the separation distance mixing time. Further
specializing to the case $\theta=1$ (corresponding to cycles of random
permutations) we prove that for $c>0$ fixed,
\[ \lim_{n \rightarrow \infty} s^*(cn^2) = 2 \sum_{i=2}^{\infty} (-1)^i (i^2-1)
e^{-ci^2}.\] Note that in contrast to the random transposition walk, there
is no cutoff.

Section \ref{other} treats other examples to which the methodology
applies. This includes Bernoulli-Laplace models, subspace walks, and a
Gibbs sampler walk on the hypercube. For most of these examples, the
spectrum is known by other methods, and separation distance results (at
least in continuous time) were described in \cite{DSa}. However the
hypercube example may be new, and in any case provides a nice illustration
of how of our method for writing $\frac{K^r(x,y)}{\pi(y)}$ as a sum of
non-negative terms allows one to determine the precise separation distance
asymptotics.

\section{Probabilistic background} \label{probback}

We will be concerned with the theory of finite Markov chains. Thus $X$
will be a finite set and $K$ a matrix indexed by $X \times X$ whose rows
sum to 1. Let $\pi$ be a distribution such that $K$ is reversible with
respect to $\pi$; this means that $\pi(x) K(x,y) = \pi(y) K(y,x)$ for all
$x,y$ and implies that $\pi$ is a stationary distribution for the Markov
chain corresponding to $K$.

Define $\langle f,g \rangle = \sum_{x \in X} f(x) g(x) \pi(x)$ for real
valued functions $f,g$ on $X$, and let $L^2(\pi)$ denote the space of such
functions. Then when $K$ is considered as an operator on $L^2(\pi)$ by
\[ Kf(x) := \sum_y K(x,y) f(y),\] it is self adjoint. Hence $K$ has an
orthonormal basis of eigenvectors $f_i(x)$ with $Kf_i(x) = \lambda_i
f_i(x)$, where both $f_i(x)$ and $\lambda_i$ are real. It is easily shown
that the eigenvalues satisfy $-1 \leq \lambda_{|X|-1} \leq \cdots \leq
\lambda_1 \leq \lambda_0=1$. If $|\lambda_1|,|\lambda_{|X|-1}|<1$, the
Markov chain is called ergodic.

\subsection{Total variation distance}

A common way to quantify the convergence rate of a Markov chain is using
total variation distance. Given probabilities $P,Q$ on $X$, one defines
the total variation distance between them as \[ ||P-Q||=\frac{1}{2}
\sum_{x \in X} |P(x)-Q(x)|.\] It is not hard to see that \[ ||P-Q|| =
\max_{A \subseteq X} |P(A)-Q(A)| .\] Let $K_x^r$ be the probability
measure given by taking $r$ steps from the starting state $x$. Researchers
in Markov chains are interested in the behavior of $||K_x^r - \pi||$.

Lemma \ref{genbound} is classical (see \cite{DH} for a proof) and relates
total variation distance to spectral properties of $K$. Note that the sum
does not include $i=0$.

\begin{lemma} \label{genbound}
\[ 4 ||K_x^r - \pi||^2 \leq \sum_{i=1}^{|X|-1} \lambda_i^{2r} |f_i(x)|^2 .\]
\end{lemma}

Lemma \ref{genbound} is remarkably effective and often leads to sharp
convergence rate results; we will apply it in Section \ref{schur}. The
main drawback with the bound in Lemma \ref{genbound} is that one rarely
knows all of the eigenvalues and eigenvectors of a Markov chain. In such
situations one typically bounds the total variation distance in terms of
$\max(|\lambda_1|,|\lambda_{|X|-1}|)$ and the results are much weaker.

\subsection{Separation distance}

Another frequently used method to quantify convergence rates of Markov
chains is to use separation distance, introduced by Aldous and Diaconis
\cite{AD1},\cite{AD2}. They define the separation distance of a Markov
chain $K$ started at $x$ as \[ s(r) = \max_{y} \left[ 1 -
\frac{K^r(x,y)}{\pi(y)} \right] \] and the maximal separation distance of
the Markov chain $K$ as \[ s^*(r) = \max_{x,y} \left[ 1 -
\frac{K^r(x,y)}{\pi(y)} \right]. \]

They show that the maximal separation distance has the nice properties:
\begin{itemize}
\item \[ \max_{x} ||K_x^r - \pi|| \leq s^*(r) \]
\item (monotonicity) $s^*(r_1) \leq s^*(r_2)$, $r_1 \geq r_2$
\item (submultiplicativity) $s^*(r_1+r_2) \leq s^*(r_1)s^*(r_2)$
\end{itemize} For every $\epsilon>0$, let $n^*_{\epsilon}$ be the smallest
number such that $s^*(n_{\epsilon}) \leq \epsilon$. Many authors consider
$n^*_{\frac{1}{2}}$ to be a definition of the separation distance mixing
time (see \cite{Pa} and references therein), and we also adopt this
convention. Heuristically, the separation distance is $\frac{1}{2}$ after
$n^*_{\frac{1}{2}}$ steps and then decreases exponentially.

Lemma \ref{tbound} will give useful upper and lower bounds for
$n^*_{\frac{1}{2}}$. It is essentially a reformulation of Corollary 2.2.9
of \cite{Pa}. By the general theory in \cite{AD2}, the random variable $T$
in Lemma \ref{tbound} always exists, but could be hard to construct.

\begin{lemma} \label{tbound} Suppose that $T$ is a random variable which takes
values in the natural numbers and satisfies $s^*(r) = \pp(T>r)$ for all $r
\geq 0$. Then
\[ \frac{\ee[T]}{2} \leq n^*_{\frac{1}{2}} \leq 2 \ee[T].\]
\end{lemma}

\begin{proof} The upper bound follows since $\pp(T>2 \ee[T]) \leq
\frac{1}{2}$. For the lower bound, note that \[ \ee[T] = \sum_{r \geq 0}
\pp(T>r) = \sum_{r \geq 0} s^*(r) \leq k + ks^*(k) + k s^*(k)^2 + \cdots =
\frac{k}{1-s^*(k)} .\] The inequality used monotonicity and
submultiplicativity. Thus if $k<\frac{\ee[T]}{2}$, then
$s^*(k)>\frac{1}{2}$, which completes the proof. \end{proof}

For the next proposition it is useful to define the distance $dist(x,y)$
between $x,y \in X$ as the smallest $r$ such that $K^r(x,y)>0$. For the
special case of birth-death chains on the set $\{0,1,\cdots,d\}$,
Proposition \ref{noeigenvec} appeared in \cite{DF} and \cite{Br}.

\begin{prop} \label{noeigenvec} (\cite{F3})
Let $K$ be a reversible ergodic Markov on a finite set $X$. Let
$1,\lambda_1,\cdots,\lambda_d$ be the distinct eigenvalues of $K$. Suppose
that $x,y$ are elements of $X$ with $dist(x,y)=d$. Then for all $r \geq
0$,
\[ 1 - \frac{K^r(x,y)}{\pi(y)} =  \sum_{i=1}^d \lambda_i^r
\left[ \prod_{j \neq i}
\frac{1-\lambda_j}{\lambda_i-\lambda_j} \right].\]
\end{prop}

The relevance of Proposition \ref{noeigenvec} to separation distance is
that one might hope that $s^*(r)$ is attained by $x,y$ satisfying
$dist(x,y)=d$. Then Proposition \ref{noeigenvec} would give an expression
for $s^*(r)$ using only the eigenvalues of $K$. Diaconis and Fill
\cite{DF} show (for $s(r)$ when the walk starts at $0$) that this hope is
realized if $K$ is a stochastically monotone birth death-chain. In the
current paper we give higher dimensional examples.

Proposition \ref{geomeig} gives a probabilistic interpretation for the
right hand side of the equation in Proposition \ref{noeigenvec}. We use
the convention that if $X$ is geometric with parameter (probability of
success) $p$, then $\pp(X=n) = p (1-p)^{n-1}$ for all $n \geq 1$.

\begin{prop} \label{geomeig} Suppose that $T=\sum_{i=1}^d X_i$,
where the random variables $X_i$ are independent, and $X_i$ is geometric
with parameter $1 - \lambda_i \in (0,1]$. If the $\lambda$'s are distinct,
then
\[ \pp(T>r) = \sum_{i=1}^d \lambda_i^r \left[ \prod_{j \neq i}
\frac{1-\lambda_j}{\lambda_i-\lambda_j} \right] \] for all natural numbers
$r$.
\end{prop}

\begin{proof} By independence, the Laplace transform of $T$ is \[
\ee[e^{-sT}] = \prod_{i=1}^d \ee[e^{-sX_i}] = \prod_{i=1}^d
\frac{1-\lambda_i}{e^s - \lambda_i}.\] Since the Laplace transform of $T$
is \[ \sum_{k \geq 1} [\pp(T>k-1)-\pp(T>k)] e^{-sk},\] it suffices to
substitute in the claimed expression for $\pp(T>k)$ and verify that one
obtains $\prod_{i=1}^d \frac{1-\lambda_i}{e^s-\lambda_i}$. Observe that
\begin{eqnarray*} & & \sum_{k \geq 1} e^{-sk} \sum_{i=1}^d \left(
\lambda_i^{k-1}- \lambda_i^{k} \right) \prod_{j \neq i}
\frac{1-\lambda_j}{\lambda_i-\lambda_j}\\ & = & \sum_{i=1}^d (1-\lambda_i)
\sum_{k \geq 1} \lambda_i^{k-1} e^{-sk} \prod_{j \neq i}
\frac{1-\lambda_j}{\lambda_i-\lambda_j} \\ & = & \sum_{i=1}^d
\frac{1-\lambda_i}{e^s-\lambda_i} \prod_{j \neq i}
\frac{1-\lambda_j}{\lambda_i-\lambda_j} \\ & = & \prod_{k=1}^d
\frac{1-\lambda_k}{e^s- \lambda_k} \sum_{i=1}^d \prod_{j \neq i}
\frac{e^s-\lambda_j}{\lambda_i-\lambda_j}. \end{eqnarray*} Letting
$t=e^s$, note that the polynomial \[ \sum_{i=1}^d \prod_{j \neq i}
\frac{e^s-\lambda_j}{\lambda_i-\lambda_j} \] is of degree at most $d-1$ in
$t$ but is equal to $1$ when $t=\lambda_i$ for $1 \leq i \leq d$. Thus the
polynomial is equal to 1, and the result follows. \end{proof}

{\it Remarks:}

\begin{enumerate}
\item Proposition \ref{geomeig} has a continuous analog where the
geometrics are exponentials \cite{BS}, and the above proof is a discrete
version of theirs.

\item For stochastically monotone birth-death chains with non-negative eigenvalues,
Proposition \ref{noeigenvec} and Proposition \ref{geomeig} lead to the
equality $s(r)=\pp(T>r)$. Here $s(r)$ is the separation distance of the
walk started at $0$, and $T$ is the sum of independent geometrics with
parameters $1-\lambda_i$, where the $\lambda_i$'s are the distinct
eigenvalues of the chain not equal to 1. This equality was first proved in
\cite{DF} using the theory of strong stationary times, and was beautifully
applied to study the cutoff phenomenon in \cite{DSa}.

\end{enumerate}

\subsection{Cut-off phenomenon}

Since the term is mentioned a few times in this article, we give a precise
definition of the cutoff phenomenon. A nice survey of the subject is
\cite{D}; we use the definition from \cite{Sal}. Consider a family of
finite sets $X_n$, each equipped with a stationary distribution $\pi_n$,
and with another probability measure $p_n$ that induces a random walk on
$X_n$. One says that there is a total variation cutoff for the family
$(X_n,\pi_n)$ if there exists a sequence $(t_n)$ of positive reals such
that
\begin{enumerate}
\item $\lim_{n \rightarrow \infty} t_n = \infty$;
\item For any $\epsilon \in (0,1)$ and $r_n = \lfloor (1+\epsilon)t_n \rfloor$,
$\lim_{n \rightarrow \infty} ||p_n^{r_n}-\pi_n||=0$;
\item For any $\epsilon \in (0,1)$ and $r_n = \lfloor (1-\epsilon)t_n \rfloor$,
$\lim_{n \rightarrow \infty} ||p_n^{r_n}-\pi_n||=1$.
\end{enumerate} For the definition of a separation cutoff, one replaces $||p_n^{
r_n}-\pi_n||$ by $s^*(r_n)$.

\section{Down-up Markov chains} \label{downup}

This section recalls the construction of down-up Markov chains on
branching diagrams and describes some main examples to be studied later in
the paper. Down-up chains appeared in \cite{F1} and more recently in
\cite{BO2}; they are obtained by composing down and up Markov chains of
Kerov \cite{K}.

The basic set-up is as follows. One starts with a branching diagram;
that is an oriented graded graph $\Gamma= \cup_{n \geq 0} \Gamma_n$
such that

\begin{enumerate}
\item $\Gamma_0$ is a single vertex $\emptyset$.
\item If the starting vertex of an edge is in $\Gamma_i$, then its end
vertex is in $\Gamma_{i+1}$.
\item Every vertex has at least one outgoing edge.
\item All $\Gamma_i$ are finite.
\end{enumerate}

For two vertices $\lambda, \Lambda \in \Gamma$, one writes $\lambda
\nearrow \Lambda$ if there is an edge from $\lambda$ to $\Lambda$.
Part of the underlying data is a multiplicity function
$\kappa(\lambda,\Lambda)$. Letting the weight of a path in $\Gamma$
be the product of the multiplicities of its edges, one defines the
dimension $d_{\Lambda}$ of a vertex $\Lambda$ to be the sum of the
weights over all maximal length paths from $\emptyset$ to $\Lambda$;
$dim(\emptyset)$ is taken to be $1$.

A set $\{M_n\}$ of probability distributions on $\Gamma_n$ is called {\it
coherent} if \[ M_{n}(\lambda) = \sum_{\Lambda: \lambda \nearrow \Lambda}
\frac{d_{\lambda} \kappa(\lambda,\Lambda)}{d_{\Lambda}}
M_{n+1}(\Lambda).\] Letting $\{M_n\}$ be a coherent set of probability
distributions, one can define the ``up'' Markov chain which transitions
from $\tau \in \Gamma_{n-1}$ to $\rho \in \Gamma_n$ with probability
$\frac{d_{\tau} M_n(\rho) \kappa(\tau,\rho)}{d_{\rho} M_{n-1}(\tau)}$.
This preserves the set $\{M_n\}$ in the sense that if $\tau$ is
distributed from $M_{n-1}$, then $\rho$ is distributed from $M_{n}$.
Similarly, one can define the ``down'' Markov chain which transitions from
$\lambda \in \Gamma_n$ to $\tau \in \Gamma_{n-1}$ with probability
$\frac{d_{\tau} \kappa(\tau,\lambda)}{d_{\lambda}}$. This also preserves
$\{M_n\}$. Composing these Markov chains by moving down and then up, one
obtains the ``down-up'' Markov chain in the level $\Gamma_n$ of the
branching diagram. This moves from $\lambda$ to $\rho$ with probability
\[ \frac{M_n(\rho)}{d_{\lambda} d_{\rho}} \sum_{\tau \in \Gamma_{n-1}}
\frac{d_{\tau}^2 \kappa(\tau, \lambda) \kappa(\tau,\rho)}
{M_{n-1}(\tau)}.\] This Markov chain has $M_n$ as its stationary
distribution and is in fact reversible with respect to $M_n$.

The reader may wonder whether there are interesting examples of coherent
probability distribution on branching diagrams. In fact there are many
such; see the surveys \cite{K} and \cite{BO1}. To make the above
definitions more concrete, we now describe two examples which are analyzed
in this paper (Young and Kingman lattices). We will also analyze down-up
walk on the Schur and Pascal lattices, but define them later.

{\it Example 1: Young lattice}

Here $\Gamma_n$ consists of all partitions of size $n$, and (identifying a
partition with its diagram in the usual way \cite{Mac}) a partition
$\lambda$ of size $n$ is adjoined to a partition $\Lambda$ of size $n+1$
if $\Lambda$ can be obtained from $\lambda$ by adding a box to some corner
of $\lambda$. The multiplicity function $\kappa(\lambda,\Lambda)$ is equal
to 1 on each edge. The dimension function $d_{\lambda}$ has an algebraic
interpretation as the dimension of the irreducible representation of the
symmetric group parameterized by $\lambda$, and there is an explicit
formula for $d_{\lambda}$ in terms of hook-lengths \cite{Sag}.

An important example of a coherent set of probability distributions on the
Young lattice is given by the so called z-measures. This is defined using
two complex parameters $z,z'$ such that $zz' \not \in \{0,-1,-2,\cdots\}$,
and assigns a partition $\lambda$ weight \[ M_{n}(\lambda) =
\frac{\prod_{b \in \lambda} (z+c(b))(z'+c(b))}{zz'(zz'+1) \cdots
(zz'+n-1)} \frac{d_{\lambda}^2}{n!}. \] Here $c(b)=j-i$ is known as the
``content'' of the box $b=(i,j)$ with row number $i$ and column number
$j$. In order that $M_n$ be strictly positive for all $n$, it is necessary
and sufficient that $(z,z')$ belongs to one of the following two sets:

\begin{itemize}
\item {\it Principal series}: Both $z,z'$ are not real and are conjugate to
each other.

\item {\it Complementary series}: Both $z,z'$ are real and are
contained in the same open interval of the form $(m,m+1)$ where $m \in
\mathbb{Z}$.
\end{itemize}

The z-measures are fundamental objects in representation theory (see
\cite{KOV1},\cite{KOV2}) and become the Plancherel measure of the
symmetric group in the limit $z,z' \rightarrow \infty$.

{\it Example 2: Kingman lattice}

Here the branching diagram is the same as the Young lattice, but the
multiplicity function $\kappa(\lambda,\Lambda)$ is the number of rows of
length $j$ in $\Lambda$, where $\lambda$ is obtained from $\Lambda$ by
removing a box from a row of length $j$. The dimension function has the
explicit form $d_{\lambda}=\frac{n!}{\lambda_1 ! \cdots \lambda_l!}$ where
$l$ is the number of rows of $\lambda$ and $\lambda_i$ is the length of
row $i$ of $\lambda$.

The Pitman distributions form a coherent set of probability distributions
on $\Gamma_n$. These are defined in terms of two parameters $\theta>0$ and
$0 \leq \alpha <1$. The Pitman distribution assigns $\lambda$ probability
\[ M_n(\lambda) = \frac{\theta(\theta+\alpha) \cdots
(\theta+(l(\lambda)-1)\alpha)}{\theta(\theta+1) \cdots (\theta+n-1)}
\frac{n!}{\prod_k m_k(\lambda)! \prod_{i=1}^{l(\lambda)} \lambda_i!}
\prod_{(i,j) \in \lambda \atop j \geq 2} (j-1-\alpha).\] Here
$m_i(\lambda)$ is the number of parts of $\lambda$ of size $i$. When
$\alpha=0$, this becomes the Ewens distribution of population genetics.
Further specializing to $\alpha=0,\theta=1$, gives that $M_n(\lambda)$ is
equal to the chance that a random permutation on $n$ symbols has cycle
type $\lambda$.

\section{Commutation relations} \label{dfrep}

It is assumed that the reader is familiar with the concept of partially
ordered sets, or posets for short. Background on posets can be found in
Chapter 3 of the text \cite{St1}. All posets considered here are assumed
to be locally finite (every interval $[x,y]$ of $P$ consists of a finite
number of elements) and graded (every maximal chain from a point $x$ to a
point $y$ has length depending only on $x,y$). It is also assumed that $P$
has an element $\hat{0}$ satisfying $x \geq \hat{0}$ for all $x \in P$.

Given a locally finite poset $P$ and $x \in P$, let $\mathbb{C}P$ denote
the complex vector space with basis $P$, and let $\mathbb{C}P_n$ denote
the subspace of $\mathbb{C}P$ spanned by the rank $n$ elements (the rank
of an element $x$ is the length $l$ of the longest chain $x_0 < x_1 <
\cdots < x_l=x$ in $P$ with top element $x$). Write $x \nearrow y$ if $y$
covers $x$ in the poset $P$. Stanley \cite{St2} defined up and down
operators $U,D$ by the condition that for $x \in P$,
\[ Ux = \sum_{y: x \nearrow y} y \ , \ Dx = \sum_{y: y \nearrow x} y.\]
These operators can be extended by linearity to $\mathbb{C}P$. For $A:
\mathbb{C}P \mapsto \mathbb{C}P$, let $A_n$ denote the restriction of $A$
to $\mathbb{C}P_n$. Notation such as $AB_n$ is unambiguous since $A(B_n)$
and $(AB)_n$ have the same meaning. Linear transformations will operate
right-to-left, e.g. $DUv=D(Uv)$, and $I$ will denote the identity
operator.

Stanley (loc. cit.) defined a locally finite, graded poset with $\hat{0}$
element to be {\it differential} if its up and down operators satisfy the
commutation relation \[ DU-UD = rI \] for some positive integer $r$. He
determined the spectrum and eigenspaces (though typically not
eigenvectors) of the operator $UD_n$. In the follow-up paper \cite{St3},
Stanley extended his ideas to the commutation relation \[ D_{n+1}U_n -
U_{n-1} D_n = r_n I_n \] where the $r_n$'s are integers.

We study the more general case that $U_n:\mathbb{C}P_n \mapsto
\mathbb{C}P_{n+1}$ and $D_n:\mathbb{C}P_n \mapsto \mathbb{C}P_{n-1}$ are
linear operators satisfying the commutation relation (1.1) of the
introduction:
\[ D_{n+1} U_n = a_n U_{n-1}D_n + b_n I_n, \] where $a_n,b_n$
are real numbers. The results we need do not all appear in \cite{Fo} (who
also studied this relation), so we briefly give statements and proofs.
This serves both to make the paper self-contained and to illustrate the
power of Stanley's methods.

Theorem \ref{eigenvalgen} determines the spectrum of $UD_n$. It can be
easily derived from Theorem 1.6.5 of \cite{Fo}.

\begin{theorem} \label{eigenvalgen} Suppose that the commutation relations
(1.1) hold and that $a_n> 0$ for all $n \geq 1$. Let $p_j$ denote the
number of elements of $P$ of rank $j$. Then the eigenvalues of $UD_n$
are $$\left\{ \begin{array}{ll} 0 & \mbox{multiplicity} \ p_n-p_{n-1} \\
\sum_{j=i}^{n-1} b_j \prod_{k=j+1}^{n-1} a_k & \mbox{multiplicity} \
p_i-p_{i-1} \ (0 \leq i \leq n-1)
\end{array} \right.$$ In particular, if $b_i=1-a_i$ for all $i$,
these become
$$\left\{ \begin{array}{ll} 0 & \mbox{multiplicity} \ p_n-p_{n-1} \\ 1 -
\prod_{k=i}^{n-1} a_k & \mbox{multiplicity} \ p_i-p_{i-1} \ (0 \leq i \leq
n-1)
\end{array} \right.$$ \end{theorem}

\begin{proof} The proof is by induction on $n$. Let $Ch(A)=Ch(A,\lambda)$
be the characteristic polynomial $det(\lambda I-A)$ of an operator $A$.
Since $Ch(U_{-1}D_0)=\lambda$, the theorem is true for $n=0$. Suppose that
$A:V \mapsto W$ and $B:W \mapsto V$ are linear transformations on finite
dimensional vector spaces $V$ and $W$ and that $dim(V)=v$ and $dim(W)=w$.
Then (by \cite{Wk}, Ch.1, Sec. 51), \[ Ch(BA)=\lambda^{v-w} Ch(AB).\]
Applying this to $D_{n+1}$ and $U_n$ gives that \begin{eqnarray*}
Ch(U_nD_{n+1},\lambda) & = & \lambda^{p_{n+1}-p_n} Ch(D_{n+1}U_n,\lambda)
\\ & = & \lambda^{p_{n+1}-p_n} Ch(a_n U_{n-1}D_n + b_n I_n, \lambda)\\ & =
& \lambda^{p_{n+1}-p_n} Ch(a_n U_{n-1}D_n, \lambda-b_n)\\ & = &
\lambda^{p_{n+1}-p_n} a_n^{p_n} det \left[ \left( \frac{\lambda-b_n}{a_n}
\right) I_n - U_{n-1}D_n \right]. \end{eqnarray*} Hence $0$ is an
eigenvalue with multiplicity at least $p_{n+1}-p_n$, and if $\lambda_k$ is
an eigenvalue of $U_{n-1}D_n$ of multiplicity $m_k$, then $a_n \lambda_k +
b_n$ is an eigenvalue of $U_nD_{n+1}$ of multiplicity at least $m_k$. This
implies the eigenvalue formula in terms of the $a,b$ variables. If one
sets $b_i=1-a_i$ for all $i$, then the sum telescopes, yielding the second
formula. \end{proof}

To compute the eigenspaces of $UD_n$, the following lemma is useful. These
eigenspaces won't be needed elsewhere in the paper, although knowing them
could prove useful in the search for eigenvectors, which by Lemma
\ref{genbound} are useful for the study of total variation distance
convergence rates.

\begin{lemma} \label{inject} Suppose that the commutation relations (1.1)
hold, with $b_0=1$ and $a_n,b_n > 0$ for all $n \geq 1$. Then the maps
$U_{n}$ are injective and the maps $D_{n+1}$ are surjective.
\end{lemma}

\begin{proof} The case $n=0$ is clear since $D_1U_0=I_0$. For $n \geq 1$, recall the
commutation relation \[ D_{n+1} U_{n} = a_n U_{n-1}D_n + b_n I_n.\] By
Theorem \ref{eigenvalgen} and the assumption that $a_n,b_n>0$ for all $n
\geq 1$, it follows that all eigenvalues of $U_{n-1}D_n$ are non-negative.
Thus all eigenvalues of $D_{n+1} U_{n}$ are positive. Thus $0$ is not an
eigenvalue and the result follows. \end{proof}

\begin{theorem} \label{eigspacevar} Suppose that the commutation relations
(1.1) hold, with $b_0=1$ and $a_n,b_n > 0$ for all $n \geq 1$. Let
$E_n(\lambda)$ denote the eigenspace of $UD_n$ corresponding to the
eigenvalue $\lambda$.
\begin{enumerate}
\item $E_n(0) = ker(D_n)$.
\item $E_n(\sum_{j=i}^{n-1} b_j \prod_{k=j+1}^{n-1} a_k) = U^{n-i}E_i(0)$.
\end{enumerate}
\end{theorem}

\begin{proof} The first assertion is clear from Lemma \ref{inject}. To
prove the second assertion, we show that \[ U_{n-1} E_{n-1} \left(
\sum_{j=i}^{n-2} b_j \prod_{k=j+1}^{n-2} a_k \right) = E_n \left(
\sum_{j=i}^{n-1} b_j \prod_{k=j+1}^{n-1} a_k \right).\] By Theorem
\ref{eigenvalgen}, the multiplicity of $\sum_{j=i}^{n-2} b_j
\prod_{k=j+1}^{n-2} a_k$ as an eigenvalue of $UD_{n-1}$ is the
multiplicity of $\sum_{j=i}^{n-1} b_j \prod_{k=j+1}^{n-1} a_k$ as an
eigenvalue of $UD_n$. Thus since $U_{n-1}$ is injective (Lemma
\ref{inject}), it is enough to check that
\[ U_{n-1} E_{n-1} \left( \sum_{j=i}^{n-2} b_j \prod_{k=j+1}^{n-2} a_k \right)
\subseteq E_n \left( \sum_{j=i}^{n-1} b_j \prod_{k=j+1}^{n-1} a_k
\right).\] So suppose that $v \in E_{n-1}(\sum_{j=i}^{n-2} b_j
\prod_{k=j+1}^{n-2} a_k)$. Then commutation relation (1.1) yields that
\begin{eqnarray*} UD_n (U_{n-1} v) & = & a_{n-1} U_{n-1} (U_{n-2} D_{n-1}
v) + b_{n-1} U_{n-1} v
\\ & = & a_{n-1} \sum_{j=i}^{n-2} b_j \prod_{k=j+1}^{n-2} a_k \cdot
U_{n-1}v + b_{n-1} U_{n-1} v
\\ & = & \sum_{j=i}^{n-1} b_j \prod_{k=j+1}^{n-1} a_k \cdot U_{n-1}v.
\end{eqnarray*}
\end{proof}

Another tool we need is an expression for $(UD)_n^r$ as a linear
combination of $(U^kD^k)_n$, extending that of \cite{St2} for the case of
differential posets.

\begin{lemma} \label{elem} Suppose that the commutation relations (1.1) hold.
Then \[ D^k U_n =  \prod_{j=n-k+1}^n a_j \cdot UD_n^k + \sum_{j=n-k+1}^n
b_j \prod_{l=j+1}^n a_l \cdot D_n^{k-1} \] for all $1 \leq k \leq n$. In
particular, if $b_i=1-a_i$ for all $i$, this becomes \[ D^k U_n =
\prod_{j=n-k+1}^n a_j \cdot UD_n^k + \left( 1 - \prod_{j=n-k+1}^n a_j
\right) \cdot D_n^{k-1}.\]
\end{lemma}

\begin{proof} This is straightforward to verify by induction on $k$,
writing $D^k=D(D^{k-1}U_n)$ and then using commutation relation (1.1).
\end{proof}

Now the desired expansion of $(UD)^r_n$ can be obtained. We remark that
for the examples studied in this paper, the coefficients in the expansion
will be non-negative.

\begin{prop} \label{pos} Suppose that the commutation relation (1.1)
holds. Then \[ (UD)_n^r = \sum_{k=0}^n A_n(r,k) (U^kD^k)_n
\] where $A_n(r,k)$ is determined by the recurrence \[ A_n(r,k) =
A_n(r-1,k-1) \prod_{j=n-k+1}^{n-1} a_j + A_n(r-1,k) \sum_{j=n-k}^{n-1} b_j
\prod_{l=j+1}^{n-1} a_l
\] with initial conditions $A_n(0,0)=1$ and $A_n(0,m)=0$ for $m \neq 0$.
In particular, if $0 \leq a_i \leq 1$, $b_i=1-a_i$ for all $i$, then the
recurrence becomes \[ A_n(r,k) = A_n(r-1,k-1) \prod_{j=n-k+1}^{n-1} a_j +
A_n(r-1,k) \left( 1 - \prod_{j=n-k}^{n-1} a_j \right) \] and all
$A_n(r,k)$ are non-negative.
\end{prop}

\begin{proof} The proposition is proved by induction on $r$. The base
case $r=0$ is clear. First applying the induction hypothesis and then
Lemma \ref{elem} yields that $(UD)_n^r$ is equal to \begin{eqnarray*}
& & \sum_{k=0}^n A_n(r-1,k) U^{k}D^{k} UD_n \\ & = & \sum_{k=0}^n A_n(r-1,k) \\
& & \cdot \left[ \prod_{j=n-k}^{n-1} a_j \cdot (U^{k+1} D^{k+1})_n +
\sum_{j=n-k}^{n-1} b_j \prod_{l=j+1}^{n-1} a_l \cdot (U^{k}D^{k})_n
\right]
\end{eqnarray*} This implies the recurrence \[ A_n(r,k) = A_n(r-1,k-1)
\prod_{j=n-k+1}^{n-1} a_j + A_n(r-1,k)  \sum_{j=n-k}^{n-1} b_j
\prod_{l=j+1}^{n-1} a_l,\] and the rest of the proposition follows
immediately.
\end{proof}

As a final result, we give a generating function for the $A_n(r,k)$ of
Proposition \ref{pos}. By comparing with Theorem \ref{eigenvalgen} one
sees that the eigenvalues of $UD_n$ appear in the generating function.

\begin{prop} \label{gen} For $k \geq 0$ set $F_k(x)=\sum_{r \geq 0}
x^r A_n(r,k)$, where $A_n(r,k)$ was defined in the statement of
Proposition \ref{pos}. Then \[ F_k(x) = \frac{x^k \prod_{i=1}^k
\prod_{j=n-i+1}^{n-1} a_j}{\prod_{i=1}^k \left( 1-x \sum_{j=n-i}^{n-1} b_j
\prod_{l=j+1}^{n-1} a_l \right)} .\] In particular, if $b_i=1-a_i$ for all
$i$, then \[ F_k(x) = \frac{x^k \prod_{i=1}^k \prod_{j=n-i+1}^{n-1}
a_j}{\prod_{i=1}^k \left[ 1-x \left( 1- \prod_{j=n-i}^{n-1} a_j \right)
\right]} .\]
\end{prop}

\begin{proof} Clearly $F_0(x)=1$. For $k \geq 1$, multiply both sides of
the recurrence of Proposition \ref{pos} by $x^r$ and sum over $r \geq 0$
to obtain that
\begin{eqnarray*}
F_k(x) & = & A_n(0,k) + \sum_{r \geq 1} x^r A_n(r,k)\\
& = & \sum_{r \geq 1} x^r A_n(r-1,k-1) \prod_{j=n-k+1}^{n-1} a_j\\ & & +
\sum_{r
\geq 1} x^r A_n(r-1,k) \sum_{j=n-k}^{n-1} b_j \prod_{l=j+1}^{n-1} a_l\\
& = & x F_{k-1}(x) \prod_{j=n-k+1}^{n-1} a_j + x F_k(x) \sum_{j=n-k}^{n-1}
b_j \prod_{l=j+1}^{n-1} a_l. \end{eqnarray*} Thus \[ F_k(x) = \frac{x
F_{k-1}(x) \prod_{j=n-k+1}^{n-1} a_j}{1 - x \sum_{j=n-k}^{n-1} b_j
\prod_{l=j+1}^{n-1} a_l},\] and the result follows by induction.
\end{proof}

\section{The Young lattice}
\label{young}

The purpose of this section is to use commutation relations to study
separation distance for down-up walk on the Young lattice. At the end of
the section, it is shown that the same asymptotics hold for up-down walk.

The setting is that of Example 1 in Section \ref{downup}. Thus the down-up
walk is on partitions of size $n$, and the chance of moving from $\lambda$
to $\rho$ is equal to \[ \frac{d_{\rho}}{n d_{\lambda}} \sum_{|\tau|=n-1
\atop \tau \nearrow \lambda,\rho}
\frac{(z+c(\rho/\tau))(z'+c(\rho/\tau))}{(zz'+n-1)}
\] and the z-measure is its stationary distribution. Here $\rho/\tau$ denotes
the box of $\rho$ not contained in $\tau$, and $c(b)=j-i$ is the
``content'' of the box $b=(i,j)$. We remind the reader that it is assumed
that either $z'=\bar{z}$ with $z \in \mathbb{C} - \mathbb{R}$, or that
$z,z'$ are real and there exists $m \in \mathbb{Z}$ such that
$m<z,z'<m+1$.

In the limiting case that $z,z' \rightarrow \infty$, the stationary
distribution becomes Plancherel measure of the symmetric group. The paper
\cite{F1} determined the eigenvalues and an orthonormal basis of
eigenvectors for down-up walk in this case. Then sharp total variation
distance convergence rates for this random walk were obtained in
\cite{F2}, and separation distance asymptotics were derived in \cite{F3}.
We give new proofs of some of these results using commutation relations,
and generalizations to the setting of z-measures.

To begin, we define operators $D_n: \mathbb{C}P_n \mapsto
\mathbb{C}P_{n-1}$ and $U_n: \mathbb{C}P_n \mapsto \mathbb{C}P_{n+1}$ as
the linear extensions of
\[ D_n(\lambda) = \sum_{\tau \nearrow \lambda} \tau \ , \ U_n(\lambda) =
\sum_{\Lambda \searrow \lambda} \frac{(z+c(\Lambda/\lambda))
(z'+c(\Lambda/\lambda))}{(zz'+n)} \Lambda.\] Note that by the hypotheses
on $z,z'$, the coefficient of any partition in $D_n(\lambda)$ or $U_n
(\lambda)$ is non-negative.

The following lemma is equivalent to Lemma 4.2 of \cite{BO2} and is
essentially due to Kerov (see \cite{O1}).

\begin{lemma} \label{comm} \[ D_{n+1}U_n = a_n U_{n-1}D_n + b_n I_n \] with
$a_n=1-\frac{1}{zz'+n}$ and $b_n=1+\frac{n}{zz'+n}$.
\end{lemma}

Let $A$ be the diagonal operator on $\mathbb{C}P$ which sends $\lambda$ to
$d_{\lambda} \cdot \lambda$. Then it is clear that the down-up walk on
Young's lattice corresponds exactly to the operator $\frac{1}{n}
(AUDA^{-1})_n$.

In Corollary \ref{eig}, $p(j)$ denotes the number of partitions of $j$. By
convention, $p(0)=1$.

\begin{cor} \label{eig} The eigenvalues of the down-up walk on the nth level of
the Young lattice are $\frac{i}{n} \left( \frac{zz'+2n-i-1}{zz'+n-1}
\right)$ $(0 \leq i \leq n)$, with multiplicity equal to
$p(n-i)-p(n-i-1)$. \end{cor}

\begin{proof} This is immediate from Theorem \ref{eigenvalgen}, Lemma \ref{comm},
and the fact that the down-up walk on Young's lattice is given by
$\frac{1}{n} (AUDA^{-1})_n$.
\end{proof}

{\it Remark:} It is not difficult to see that $p(n-i)-p(n-i-1)$ is equal
to the number of partitions of $n$ with $i$ 1's. Indeed, using the
notation that $[u^n]f(u)$ is the coefficient of $u^n$ in $f(u)$, one has
that
\begin{eqnarray*} p(n-i)-p(n-i-1) & = & [u^{n-i}] \prod_{j \geq 1}
(1-u^j)^{-1} - [u^{n-i-1}] \prod_{j \geq 1}
(1-u^j)^{-1}\\ & = & [u^{n-i}] (1-u) \prod_{j \geq 1} (1-u^j)^{-1} \\
& = & [u^{n-i}] \prod_{j \geq 2} (1-u^j)^{-1}, \end{eqnarray*} which
is the number of partitions of $n-i$ with no 1's.

\begin{center}
\end{center}

Proposition \ref{sepyoung} is crucial for determining where the maximal
separation distance of down-up walk on Young's lattice is attained. Its
statement uses the notation that if $B:\mathbb{C}P \mapsto \mathbb{C}P$,
then $B[\mu,\lambda]$ is the coefficient of $\lambda$ in $B(\mu)$.

\begin{prop} \label{sepyoung} Let $\pi(\lambda)$ be the z-measure evaluated at
$\lambda$, and let $r$ be a non-negative integer. Then the quantity \[
\frac{(\frac{1}{n} AUDA^{-1})^r[\mu,\lambda]} {\pi(\lambda)} \] is
minimized (among pairs of partitions of size $n$) by
$\mu=(n),\lambda=(1^n)$ or $\mu=(1^n),\lambda=(n)$.
\end{prop}

\begin{proof} Lemma \ref{comm} and Proposition \ref{pos} give that \[
\frac{(UD)^r[\mu,\lambda]}{\pi(\lambda)} = \sum_{k=0}^n A_n(r,k)
\frac{U^kD^k[\mu,\lambda]}{\pi(\lambda)}, \] where $A_n(r,k)$ is
determined by the recursion of Proposition \ref{pos}. Thus \[
\frac{(\frac{1}{n} AUDA^{-1})^r[\mu,\lambda]}{\pi(\lambda)} =
\frac{1}{n^r} \sum_{k=0}^n \frac{d_{\lambda} A_n(r,k)
U^kD^k[\mu,\lambda]}{d_{\mu} \pi(\lambda)}.\]

The proposition follows immediately from three claims:

\begin{itemize}
\item All terms in the sum are non-negative. Indeed, since
$b_n \geq 0$ for $n \geq 0$ and $a_n \geq 0$ for $n \geq 1$, the recursion
for $A_n(r,k)$ implies that $A_n(r,k) \geq 0$. Noting that $U,D$ map
non-negative linear combinations of partitions to non-negative linear
combinations of partitions, the claim follows.

\item If $\mu=(n),\lambda=(1^n)$ or $\mu=(1^n),\lambda=(n)$, then the summands
for $0 \leq k \leq n-2$ vanish. Indeed, for such $k$ it is impossible to
move from the partition $\mu$ to the partition $\lambda$ by removing $k$
boxes one at a time and then reattaching $k$ boxes one at a time.

\item The $k=n-1$ and $k=n$ summands are independent of both $\mu$ and $\lambda$.
Indeed, for the $k=n-1$ summand one has that
\begin{eqnarray*}
 & & \frac{d_{\lambda} A_n(r,n-1) U^{n-1}D^{n-1}[\mu,\lambda]}{n^r d_{\mu} \pi(\lambda)}\\ &
= &
\frac{d_{\lambda} A_n(r,n-1) U^{n-1}[(1),\lambda]}{n^r \pi(\lambda)} \\
& = & \frac{d_{\lambda}^2 A_n(r,n-1) \prod_{b \in \lambda \atop b \neq
(1,1)} (z+c(b))(z'+c(b))}{n^r (zz'+1) \cdots (zz'+(n-1)) \pi(\lambda)} \\
& = & \frac{n! A_n(r,n-1)}{n^r}. \end{eqnarray*} The first equality used
the fact that there are $d_{\mu}$ ways to go from $\mu$ to $(1)$ by
removing a box at a time. The second equality used the fact that all
$d_{\lambda}$ ways of transitioning from $(1)$ to $\lambda$ in $n-1$
upward steps give the same contribution to $U^{n-1}[(1),\lambda]$.

A similar argument shows that the $k=n$ summand is equal to $\frac{n!
A_n(r,n)}{n^r}$.

\end{itemize}
\end{proof}

Corollary \ref{sepdistz} gives an expression for maximal separation
distance.

\begin{cor} \label{sepdistz} Let $s^*(r)$ be the maximal separation distance
after r iterations of the down-up chain $K$ on the nth level of the Young
lattice. Then $s^*(r) = \pp(T>r)$, where $T$ is a sum of independent
geometrics with parameters $1-\frac{i}{n} \left(
\frac{zz'+2n-i-1}{zz'+n-1} \right)$ for $0 \leq i \leq n-2$. \end{cor}

\begin{proof} By Proposition \ref{sepyoung}, $s^*(r) = 1 - \frac{K^r((n),(1^n))}
{\pi(1^n)}$. By Proposition \ref{eig}, the down-up walk has $n$ distinct
eigenvalues, namely $1$ and $\frac{i}{n} \left(\frac{zz'+2n-i-1}{zz'+n-1}
\right)$ for $0 \leq i \leq n-2$. Since the distance between $(n)$ and
$(1^n)$ is $n-1$, the result follows from Propositions \ref{noeigenvec}
and \ref{geomeig}.
\end{proof}

Theorem \ref{exactplan} gives a precise expression for the asymptotics of
separation distance in the special case that $z,z' \rightarrow \infty$.
Then the stationary distribution is Plancherel measure of the symmetric
group, and these asymptotics were obtained earlier in \cite{F3}. Here we
present a new proof which involves determining the numbers $A_n(r,k)$.
This technique is likely to prove useful for other problems; in
particular, we apply it again later in this paper (Proposition
\ref{cutpas}).

\begin{theorem} \label{exactplan} Let $s^*(r)$ be the maximal separation distance
after $r$ iterations of the down-up walk on the nth level of the Young
lattice, in the special case that $z,z' \rightarrow \infty$.
\begin{enumerate}
\item $s^*(r) = 1 - \frac{n! S(r,n-1)}{n^r} - \frac{n! S(r,n)}{n^r}$, where
$S(r,k)$ is a Stirling number of the second kind (i.e. the number of
partitions of an r-set into k blocks).

\item For c fixed in $\mathbb{R}$ and $n \rightarrow \infty$,
\[ s^*(n \log(n)+cn) = 1 - e^{-e^{-c}}(1+e^{-c}) + O \left( \frac{\log(n)}{n}
\right) .\]
\end{enumerate}
\end{theorem}

\begin{proof} For the first assertion, the proof of Proposition \ref{sepyoung}
gives that \[ s^*(r) = 1 - \frac{n! A_n(r,n-1)}{n^r} - \frac{n!
A_n(r,n)}{n^r}.\] The recurrence in Proposition \ref{pos} is \[ A_n(r,k) =
A_n(r-1,k-1) + k A_n(r-1,k) \] with initial conditions $A_n(0,m) =
\delta_{0,m}$. The solution to this recurrence is $A_n(r,k)=S(r,k)$ (see
also Proposition 4.9 of \cite{St2}), and the first assertion follows.

Let $P(n,r,k)$ denote the probability that when $r$ balls are dropped
uniformly at random into $n$ boxes, there are $k$ occupied boxes. It is
straightforward to see that $P(n,r,k) = \frac{S(r,k) k! {n \choose
k}}{n^r}$. Indeed, occupying $k$ boxes using $r$ balls is equivalent to
forming an ordered set partition of $\{1,\cdots,r\}$ into $k$ blocks and
then choosing $k$ of the $n$ boxes. Thus, \[ s^*(r) = 1 - P(n,r,n-1) -
P(n,r,n).\] Now we use asymptotics of the coupon collector's problem: it
follows from Section 6 of \cite{CDM} that when $n \log(n)+cn$ balls are
dropped into $n$ boxes, the number of unoccupied boxes converges to a
Poisson distribution with mean $e^{-c}$, and that the error term in total
variation distance is $O(\frac{\log(n)}{n})$. The chance that a Poisson
random variable with mean $e^{-c}$ takes value not equal to 0 or 1 is
$1-e^{-e^{-c}}(1+e^{-c})$, which completes the proof.
\end{proof}

For general values of $z,z'$ it is not evident how to obtain results as
clean as Theorem \ref{exactplan}. However Proposition \ref{someinfo} gives
explicit upper and lower bounds for the separation distance mixing time.
For $z,z'$ fixed and $n$ growing, these are both order $n^2$.

\begin{prop} \label{someinfo} Let $n^*_{\frac{1}{2}}$ be the separation
distance mixing time of down-up walk (corresponding to z-measure) on the
nth level of the Young lattice. Then $\frac{\ee[T]}{2} \leq
n^*_{\frac{1}{2}} \leq 2 \ee[T]$, where $T$ is as in Corollary
\ref{sepdistz}. Moreover, if $zz'=1$ then \[ \ee[T] = \sum_{i=2}^n
\frac{n^2}{i^2} \sim n^2 \left( \frac{\pi^2}{6} -1 \right),
\] and if $zz' \neq 1$ then \begin{eqnarray*} 1 + \frac{n(zz'+n-1)}{1-zz'} \log \left(
\frac{2(n+zz'-1)}{n(zz'+1)} \right)  \leq  \ee[T] \\ \leq
\frac{n(zz'+n-1)}{1-zz'} \log \left( \frac{n+zz'-1}{n(zz')} \right).
\end{eqnarray*}
\end{prop}

\begin{proof} By Lemma \ref{tbound}, $\frac{\ee[T]}{2} \leq
n^*_{\frac{1}{2}} \leq 2 \ee[T]$.  Linearity of expectation gives that \[
\ee[T] = \sum_{i=0}^{n-2} \frac{1}{1 - \frac{i}{n} \left(
\frac{zz'+2n-i-1}{zz'+n-1} \right)}.\] When $zz'=1$, \[ \ee[T] =
\sum_{i=0}^{n-2} \frac{1}{1 - \frac{i(2n-i)}{n^2}} = n^2 \sum_{i=0}^{n-2}
\frac{1}{(n-i)^2} = n^2 \sum_{i=2}^n \frac{1}{i^2} \sim n^2 \left(
\frac{\pi^2}{6} -1 \right).\]

For $zz' \neq 1$, the fact that $\frac{i(zz'+2n-i-1)}{n(zz'+n-1)}$ is
monotone increasing for $i \in [0,n-1]$ gives that \[ 1 + \int_{0}^{n-2}
\frac{1}{1 - \frac{t}{n} \left( \frac{zz'+2n-t-1}{zz'+n-1} \right)} dt
\leq \ee[T] \leq \int_{0}^{n-1} \frac{1}{1 - \frac{t}{n} \left(
\frac{zz'+2n-t-1}{zz'+n-1} \right)} dt.\] Consider the upper bound on
$\ee[T]$. Since $zz' \neq 1$, it is equal to \begin{eqnarray*} & &
\frac{n(zz'+n-1)}{1-zz'} \int_{0}^{n-1} \left( \frac{1}{t-n} -
\frac{1}{t-(n+zz'-1)} \right) dt \\ & = & \frac{n(zz'+n-1)}{1-zz'} \log
\left( \frac{n+zz'-1}{n(zz')} \right). \end{eqnarray*} Similarly, since
$zz' \neq 1$, the lower bound on $\ee[T]$ is equal to \begin{eqnarray*} &
& 1 + \frac{n(zz'+n-1)}{1-zz'} \int_{0}^{n-2} \left( \frac{1}{t-n} -
\frac{1}{t-(n+zz'-1)} \right) dt \\ & = & 1 + \frac{n(zz'+n-1)}{1-zz'}
\log \left( \frac{2(n+zz'-1)}{n(zz'+1)} \right) .\end{eqnarray*}
\end{proof}

To close this section we prove Proposition \ref{samerate}. It implies that
the up-down and down-up walks have the same convergence rate asymptotics.

\begin{prop} \label{samerate} Let $s^*_{UD_n}(r)$ be the maximal separation
distance after $r$ iterations of the down-up chain (corresponding to
z-measure) on the Young lattice, and let $s^*_{DU_n}(r)$ be the
corresponding quantity for the up-down chain. Then \[ s^*_{DU_n}(r) =
s^*_{UD_{n+1}}(r+1) \] for all $n \geq 1, r \geq 0$.
\end{prop}

\begin{proof} Using the notation of Proposition \ref{sepyoung}, the up-down
chain corresponds to the operator $\frac{1}{n+1} ADU_{n} A^{-1}$. Lemma
\ref{comm} implies that \[ \frac{ (\frac{1}{n+1} ADU_n A^{-1})^r
[\mu,\lambda]} {\pi(\lambda)} = \frac{1}{(n+1)^r} \sum_{l=0}^r {r \choose
l} a_n^l b_n^{r-l} n^l \frac{ (\frac{1}{n} AUD_n A^{-1})^l [\mu,\lambda]}
{\pi(\lambda)}, \] where $a_n=1 - \frac{1}{zz'+n}$ and $b_n=1 +
\frac{n}{zz'+n}$. Hence Proposition \ref{sepyoung} gives that this
quantity is minimized by $\mu=(n), \lambda=(1^n)$ or $\mu=(1^n),
\lambda=(n)$.

By Corollary \ref{eig} and Lemma \ref{comm}, the distinct eigenvalues of
the up-down chain are $1$ and $t_j:=\frac{j(zz'+2n-j+1)}{(n+1)(zz'+n)}$
where $1 \leq j \leq n-1$. Hence Proposition \ref{noeigenvec} gives that
\[ s^*_{DU_n}(r) = \sum_{j=1}^{n-1} (t_j)^r \prod_{k \neq j \atop 1 \leq k \leq n-1}
\left( \frac{1-t_k}{t_j-t_k} \right).\] On the other hand, applying
Proposition \ref{noeigenvec} to the down-up chain gives that
\[ s^*_{UD_{n+1}}(r+1) = \sum_{j=0}^{n-1} (t_j)^{r+1} \prod_{k
\neq j \atop 0 \leq k \leq n-1} \left( \frac{1-t_k}{t_j-t_k} \right).\]
Since $t_0=0$, this becomes \[ \sum_{j=1}^{n-1} (t_j)^{r+1} \prod_{k \neq
j \atop 0 \leq k \leq n-1} \left( \frac{1-t_k}{t_j-t_k} \right) =
\sum_{j=1}^{n-1} (t_j)^{r} \prod_{k \neq j \atop 1 \leq k \leq n-1} \left(
\frac{1-t_k}{t_j-t_k} \right),
\] as desired.
\end{proof}

\section{The Schur lattice} \label{schur}

In this example the underlying lattice is the Schur lattice. This is the
sublattice of Young's lattice consisting of the partitions of $n$ into
distinct parts. We show that commutation relations can be used to compute
the spectrum of down-up walk on the Schur lattice, but our approach does
not determine the separation distance convergence rate (the obstacles are
described in the second remark after Proposition \ref{eig2}). We do
however give a complete diagonalization of the Markov chain, and use it to
study the total variation distance convergence rate. The upper bound
derived here is in fact quite sharp and there is a cutoff at $\frac{1}{2}n
\log(n)$. We omit the rather involved proof of a matching lower bound but
give a careful statement and explain the proof technique in the remarks
after Theorem \ref{updown}.

It will be convenient to let $DP(n)$ denote the set of partitions of $n$
into distinct parts and $OP(n)$ denote the set of partitions of $n$ into
odd parts. Using the terminology of Section \ref{downup}, there is a
coherent set of probability distributions on the Schur lattice called the
shifted Plancherel measures. The $n$th measure chooses a partition
$\lambda \in DP(n)$ with probability \[ \pi(\lambda):=
\frac{2^{n-l(\lambda)} g_{\lambda}^2}{n!},\] where $l(\lambda)$ is the
number of parts of $\lambda$ and $g_{\lambda}$ is the number of standard
shifted tableaux of shape $\lambda$ (\cite{HH},\cite{Mac}). This measure
is of interest to researchers in asymptotic combinatorics and
representation theory \cite{B},\cite{Mat},\cite{TW}, \cite{I}.

In the terminology of Section \ref{downup}, it is known (see for instance
\cite{BO1}) that the dimension of $\lambda \in DP(n)$ is equal to
$g_{\lambda}$. Hence the down-up chain on the set $DP(n)$ transitions from
$\lambda$ to $\rho$ with probability
\[ \frac{2 g_{\rho}}{n g_{\lambda}} \sum_{\tau
\nearrow \lambda,\rho} 2^{l(\tau)-l(\rho)}.\] An application of this
Markov chain appears in \cite{F4}. However nothing seems to be known about
its convergence rate.

We will diagonalize this chain (determining eigenvalues and eigenvectors).
Before doing this we note that commutation relations can also be used to
derive its eigenvalues. The key is the following observation of Stanley
\cite {St3}. He defined down and up operators $D,U$ for the Schur lattice
by:
\[ D(\lambda) = \sum_{\mu \nearrow \lambda} \mu , \ U(\lambda) = 2
\sum_{\mu \searrow \lambda \atop l(\mu)=l(\lambda)} \mu + \sum_{\nu
\searrow \lambda \atop l(\nu)>l(\lambda)} \nu, \] and showed that they
satisfy the commutation relation \num \begin{equation} \label{stan}
D_{n+1}U_n = U_{n-1} D_n + I_n
\end{equation} for all $n \geq 0$.

In Proposition \ref{eig2}, $p^*(j)$ denotes the number of partitions of
$j$ into distinct parts.

\begin{prop} \label{eig2} The eigenvalues of the down-up walk on the Schur
lattice are $\frac{i}{n}$ $(0 \leq i \leq n)$, with multiplicity equal to
$p^*(n-i)-p^*(n-i-1)$.
\end{prop}

\begin{proof} Let $A$ be the diagonal operator on $\mathbb{C}P$ which sends
$\lambda$ to $g_{\lambda} \cdot \lambda$. It is easily seen that the
down-up chain is equivalent to the operator $\frac{1}{n} (AUDA^{-1})_n$.
The result now follows from commutation relation (\ref{stan}) and Theorem
\ref{eigenvalgen}. \end{proof}

{\it Remarks:}
\begin{enumerate}
\item It is well known that $|DP(n)|=|OP(n)|$. Using generating functions
as in the remark after Proposition \ref{eig}, one can show that
$p^*(n-i)-p^*(n-i-1)$ is equal to the number of odd partitions of $n$ with
$i$ parts equal to 1. This also follows by comparing Proposition
\ref{eig2} with Proposition \ref{diag} below.

\item From the previous remark, it is easily seen that the number of
distinct eigenvalues of $UD_n$ is $n-2$ for large enough $n$ (an odd
partition of $n$ can't have $i$ parts of size $1$ for $i=n-1,n-2,n-4$).
However the diameter of down-up walk on the Schur lattice can be smaller
than $n-3$ (for $n=8$ it is 4). This blocks the use of Proposition
\ref{noeigenvec} and also complicates the analysis of where the maximal
separation distance is attained, as the proof of Proposition
\ref{sepyoung} does not carry over.

\end{enumerate}

To upper bound the total variation distance convergence rate, the
following diagonalization of the down-up walk is crucial. The eigenvectors
are given in terms of symmetric functions, more precisely in terms of
$X^{\lambda}_{\mu}$ which is defined as the coefficient of the
Hall-Littlewood polynomial $P_{\lambda}(x;-1)$ in the power sum symmetric
function $p_{\mu}(x)$. The reader unfamiliar with these concepts can
either consult Chapter 3 of \cite{Mac} (which calls these coefficients
$X^{\lambda}_{\mu}(-1)$), or can just proceed to Theorem \ref{updown}. We
also use the notation that $z_{\mu}=\prod_i i^{m_i(\mu)} m_i(\mu)!$, where
$m_i(\mu)$ is the number of parts of $\mu$ of size $i$. This is the number
of permutations which commute with a fixed permutation of cycle type
$\mu$.

\begin{prop} \label{diag}
\begin{enumerate}
\item The eigenvalues of down-up walk on the Schur lattice are
parameterized by $\mu \in OP(n)$ and are $\frac{m_1(\mu)}{n}$, where
$m_1(\mu)$ is the number of parts of $\mu$ of size 1.
\item The functions
$\psi_{\mu}(\lambda) = \sqrt{\frac{n!}{z_{\mu} 2^{n-l(\mu)}}}
\frac{X^{\lambda}_{\mu}}{g_{\lambda}}$ are a corresponding basis of
eigenvectors, orthonormal with respect to the inner product \[ \langle
f_1,f_2 \rangle = \sum_{\lambda \in DP(n)} f_1(\lambda)
\overline{f_2(\lambda)} \frac{2^{n-l(\lambda)} g_{\lambda}^2}{n!}.\]
\end{enumerate}
\end{prop}

\begin{proof} It follows from Lemma 5.6 and Corollary 5.11 of \cite{F4}
that the $\psi_{\mu}$ are an orthonormal basis of eigenvectors with
eigenvalue $\frac{m_1(\mu)-2}{n-2}$ for a certain operator $J_{(n-1,1)}$,
defined by \[ J_{(n-1,1)}(\lambda,\rho) = \frac{g_{\rho}}{2^{l(\rho)}
g_{\lambda} (n-2)} \sum_{\nu \in OP(n)} \frac{2^{l(\nu)} X^{\lambda}_{\nu}
X^{\rho}_{\nu} (m_1(\nu)-2)}{z_{\nu}}.\] The proposition follows from the
claim that the chance that the down-up chain moves from $\lambda$ to
$\rho$ is equal to
\[ \frac{(n-2)J_{(n-1,1)}(\lambda,\rho)}{n} + \frac{2}{n}
\delta_{\lambda,\rho}
\] where $\delta_{\lambda,\rho}$ is $1$ if $\lambda=\rho$ and vanishes
otherwise. For the case that $\lambda \neq \rho$, the claim follows from
the statement of Proposition 5.9 of \cite{F4}, and for the case
$\lambda=\rho$, it follows from the proof of Proposition 5.9 and Lemma 5.3
of \cite{F4}. \end{proof}

Finally, we use the diagonalization to study total variation distance for
down-up walk on the Schur lattice.

\begin{theorem} \label{updown} Let $K^r$ denote the distribution of the
down-up walk on the Schur lattice started from $(n)$ after $r$ steps, and
let $\pi$ denote the shifted Plancherel measure. For $r=\frac{1}{2}n
\log(n)+cn$ with $c>0$, \[ ||K^r - \pi|| \leq \frac{e^{-3c}}{4}.\]
\end{theorem}

\begin{proof} The diagonalization of the down-up walk, together with Lemma
\ref{genbound} and the facts \cite{Mac} that $g_{(n)}=1$ and
$X^{(n)}_{\mu}=1$ for all $\mu$, gives that
\begin{eqnarray*}
||K^r - \pi||^2 & \leq & \frac{1}{4} \sum_{\mu \neq (1^n) \atop \mu \in
OP(n)}
\left( \frac{m_1(\mu)}{n} \right)^{2r} \frac{n!}{z_{\mu} 2^{n-l(\mu)}} \\
& = & \frac{1}{4} \sum_{i=1}^{n-2} \left( \frac{i}{n} \right)^{2r}
\sum_{\mu \in OP(n) \atop m_1(\mu)=i} \frac{n!}{z_{\mu} 2^{n-l(\mu)}}.
\end{eqnarray*} Letting $[u^n]f(u)$ denote the coefficient of $u^n$ in $f(u)$,
the cycle index of the symmetric group (reviewed in Chapter 4 of
\cite{Wi}) yields that
\begin{eqnarray*} \sum_{\mu \in OP(n) \atop m_1(\mu)=i} \frac{n!}{z_{\mu}
2^{n-l(\mu)}} & = & \frac{n!}{i! 2^{n-i}} [u^{n-i}] \prod_{m \geq 3 \atop
odd} e^{\frac{2u^m}{m}} \\ & = & \frac{n!}{i! 2^{n-i}} [u^{n-i}]
\frac{1}{e^{2u}} \prod_{m \geq 1 \atop odd} e^{\frac{2u^m}{m}} \\ & = &
\frac{n!}{i! 2^{n-i}} [u^{n-i}] \frac{(1+u)}{(1-u) e^{2u}} \\ & = &
\frac{n!}{i! 2^{n-i}} \left[ \sum_{j=0}^{n-i} \frac{(-2)^j}{j!} +
\sum_{j=0}^{n-i-1} \frac{(-2)^j}{j!} \right]. \end{eqnarray*} It is easily
checked that \[ \sum_{j=0}^{n-i} \frac{(-2)^j}{j!} + \sum_{j=0}^{n-i-1}
\frac{(-2)^j}{j!} \] vanishes if $n-i=1,2,4$ and when $n-i>0$ is at most
$2/3$.

Thus \begin{eqnarray*} ||K^r - \pi||^2 & \leq & \frac{1}{6} \sum_{i=1
\atop i \neq n-4}^{n-3} \left( \frac{i}{n} \right)^{2r} \frac{n!}{i!
2^{n-i}} \\ & \leq & \frac{1}{6} \sum_{j=3}^n (1-\frac{j}{n})^{2r}
\frac{n!}{(n-j)! 2^j}\\ & = & \frac{1}{6} \sum_{j=3}^n \frac{n!}{(n-j)!
2^j} e^{2r \cdot \log(1-j/n)} \\ & \leq & \frac{1}{6} \sum_{j=3}^n
\frac{n!}{(n-j)!} \frac{e^{-2rj/n}}{2^j} \\ & = & \frac{1}{6} \sum_{j=3}^n
\frac{n!}{(n-j)!} \frac{e^{-2cj}}{n^j 2^j} \\ & \leq & \frac{1}{6}
\sum_{j=3}^n \frac{e^{-2cj}}{2^j} \\ & = & \frac{e^{-6c}}{48(1-e^{-2c}/2)}
\\ & \leq & \frac{e^{-6c}}{24}. \end{eqnarray*} Taking square roots
completes the proof. \end{proof}

{\it Remarks:}

\begin{enumerate}
\item One can prove that there are positive universal constants $A,B$ such
that for all $c>0$ and $r=\frac{1}{2} n \log(n)-cn$ with $n$ large enough
(depending on $c$),
\[ ||K^r - \pi|| \geq 1-Ae^{-Bc}.\] The proof method is analogous to that
used in \cite{F2} for the case of Plancherel measure of the symmetric
group, but the combinatorics is more tedious. One can compute the mean and
variance of the eigenfunction $\psi_{(3,1^{n-3})}$ under both $\pi$ and
the measure $K^r$, and then deduce the lower bound from Chebyshev's
inequality.

\item From commutation relation (\ref{stan}), the results in this section
give (in the notation of Proposition \ref{diag}) that up-down walk on the
Schur lattice has eigenvalues $\frac{m_1(\mu)+1}{n+1}$ and the same
eigenfunctions as down-up walk. Arguing as in Theorem \ref{updown} gives
that the walks have the same convergence rate asymptotics.

\end{enumerate}

To conclude this section, we mention that the techniques in it can be used
to analyze total variation distance convergence rates for down-up walk on
the Jack lattice. Here the stationary distribution is the so-called
Jack$_{\alpha}$ measure on partitions, which in the special case
$\alpha=1$ gives the Plancherel measure of the symmetric group. The
importance of Jack$_{\alpha}$ measure is discussed in Okounkov \cite{O2},
and some results about it appear in \cite{BO4} and \cite{F5}. In
particular, Proposition 6.2 of \cite{F5} explicitly diagonalizes down-up
walk on the Jack lattice. The eigenvalues turn out to be independent of
$\alpha$ and are $1$ and $\frac{i}{n}$ for $0 \leq i \leq n-2$. The
eigenvectors are the coefficients of power sum symmetric functions in the
Jack polynomials with parameter $\alpha$. Further details may appear
elsewhere.

\section{The Kingman lattice} \label{kingman}

This section uses commutation relations to study down-up walk on the
Kingman lattice. The stationary distribution is the Pitman distribution
with parameters $\theta,\alpha$ where $\theta>0$ and $0 \leq \alpha <1$
(Example 2 in Section \ref{downup}). We show that the eigenvalues and
separation distance do not depend on $\alpha$ and prove order $n^2$ upper
and lower mixing time bounds. Very precise convergence rate results are
given when $\theta=1$. This is probably the most interesting case, since
when $\alpha=0, \theta=1$ the stationary distribution corresponds to the
cycle structure of random permutations.

The down-up walk studied in this section is more ``local'' the the random
transposition walk, in the sense that the underlying partition is changed
by removing a single box and then reattaching it somewhere. In the random
transposition walk, the change is more violent: two cycles can merge into
one cycle or a single cycle can be broken into two cycles. Local walks
tend to be more useful for Stein's method than non-local walks (see
\cite{R} for some rigorous results in this direction), and this down-up
walk was described in Section 2 of \cite{F1} in the context of Stein's
method. The recent paper \cite{Pe} applies down-up walk on Kingman's
lattice to define a new family of infinite dimensional diffusions, which
includes the infinitely-many-neutral-alleles-diffusion model of Ethier and
Kurtz.

Now we begin the analysis of the down-up chain corresponding to the Pitman
distribution with parameters $\theta>0$ and $0 \leq \alpha<1$. By the
formulas in Section \ref{downup}, one sees that the down chain removes a
box from a row of length $j$ with probability $\frac{j m_j(\lambda)}{n}$
and that the up chain adds a box to a row of $\lambda$ of length $k \geq
1$ with probability $\frac{(k-\alpha) m_k(\lambda)}{\theta+n}$ or to a row
of length 0 with probability $\frac{\theta+\alpha l(\lambda)}{\theta+n}$,
where $l(\lambda)$ is the number of parts of $\lambda$. In the biological
context $(\alpha=0)$, the rows of $\lambda$ could represent the count of
individuals of each type in a population. Then the down move corresponds
to the death of a random individual, and the up move corresponds to a
birth (which is the same type as the random parent or a new type with
probability $\frac{\theta}{\theta+n}$).

Let $P$ be the poset of partitions with the same partial order as in
Kingman's lattice, where we disregard edge multiplicities; this is the
same partial order as in Young's lattice. It is natural to define
operators $D,U: \mathbb{C}P \mapsto \mathbb{C}P$ as follows. The
coefficient of $\tau$ in $D_n(\lambda)$ is defined to be the probability
that from $\lambda$, the down-chain transitions to $\tau$. The coefficient
of $\Lambda$ in $U_n(\lambda)$ is defined to be the probability that from
$\lambda$, the up-chain transitions to $\Lambda$. Thus the down-up walk on
Kingman's lattice arising from Pitman's distribution is just the operator
$UD_n$.

The following commutation relation is crucial. Note that a closely related
commutation relation appears in \cite{Pe}.

\begin{prop} \label{commutking} Consider down-up walk on the Kingman lattice
with parameters $\theta>0$ and $0 \leq \alpha <1$. Letting $a_n =
\frac{n(\theta+n-1)}{(n+1)(\theta+n)}$, one has that
\[ D_{n+1} U_n = a_n U_{n-1} D_n + (1-a_n) I_n, \] for all $n \geq 0$.
\end{prop}

\begin{proof} First we consider the case that $\lambda,\rho$ are distinct
partitions of $n$. Then in order to move from $\lambda$ to $\rho$ by going
up and then going down, one must add a box to a row of length $k$ of
$\lambda$ and then remove a box from a row of length $j$. Similarly, in
order to move from $\lambda$ to $\rho$ by going down and then going up,
one must remove a box from a row of length $j$ of $\lambda$, and then add
a box to a row of length $k$. In both situations one has that $j \geq 1, k
\geq 0$ and $j \neq k+1$. From this it is straightforward to check
(treating separately the cases that $k>0$ and $k=0$), that the coefficient
of $\rho$ in \[ (n+1)(\theta+n) DU_n(\lambda) - n(\theta+n-1)
UD_n(\lambda)
\] is $0$.

The second case to consider is that $\lambda = \rho$ are the same
partition of $n$. Then $j=k+1$, and the coefficient of $\lambda$ in
$(n+1)(\theta+n) DU_n(\lambda)$ is \[ [\theta+\alpha l(\lambda)]
[m_1(\lambda)+1] + \sum_{k \geq 1} [(k-\alpha) m_k(\lambda)]
[(k+1)(m_{k+1}(\lambda)+1)].\] Similarly, the coefficient of $\lambda$ in
$n(\theta+n-1)UD_n(\lambda)$ is \[
m_1(\lambda)[\theta+\alpha(l(\lambda)-1)] + \sum_{k \geq 1}
[(k+1)m_{k+1}(\lambda)] [(k-\alpha)(m_k(\lambda)+1)] .\] Hence the
coefficient of $\lambda$ in \[ (n+1)(\theta+n) D U_n(\lambda) -
n(\theta+n-1) U D_n(\lambda) \] is
\begin{eqnarray*} & & \theta + \alpha l(\lambda) + \alpha m_1(\lambda) +
 \sum_{k \geq 1} (k-\alpha)(k+1)
(m_k(\lambda)-m_{k+1}(\lambda))\\ & = & \theta + \alpha l(\lambda) -
\alpha m_1(\lambda) + 2m_1(\lambda)\\ & & + \sum_{k \geq 2}
m_k(\lambda)[(k+1)(k-\alpha)-k(k-1-\alpha)]\\
& = & \theta + \alpha l(\lambda) - \alpha m_1(\lambda) + 2m_1(\lambda) +
\sum_{k \geq 2} (2k-\alpha) m_k(\lambda) \\
& = & \theta + 2n. \end{eqnarray*}
\end{proof}

Corollary \ref{eigenking} determines the eigenvalues of the down-up walk
on the Kingman lattice with parameters $\theta,\alpha$. It is interesting
that these are independent of the parameter $\alpha$. We remark that since
$p(1)=p(0)=1$, the eigenvalue $1-\frac{\theta}{n(\theta+n-1)}$ in
Corollary \ref{eigenking} has multiplicity 0.

\begin{cor} \label{eigenking} Let $p(j)$ denote the number of integer
partitions of $j$. Then the eigenvalues of $U D_n$ are $1 -
\frac{i(\theta+i-1)}{n(\theta+n-1)}$ with multiplicity $p(i)-p(i-1)$ $(0
\leq i \leq n)$.
\end{cor}

\begin{proof} This is immediate from Theorem \ref{eigenvalgen} and
Proposition \ref{commutking}. \end{proof}

Next we will study maximal separation distance for the down-up walk on the
Kingman lattice. The first step is to determine where this is attained.
Given a linear operator $B:\mathbb{C}P \mapsto \mathbb{C}P$, and
partitions $\mu,\lambda$, it is convenient to let $B[\mu,\lambda]$ denote
the coefficient of $\lambda$ in $B(\mu)$.

\begin{prop} \label{sepking} Let $\pi$ be the Pitman distribution with
parameters $\theta>0$ and $0 \leq \alpha <1$. Let $r$ be a non-negative
integer. The quantity $\frac{(UD)^r [\mu,\lambda]}{\pi(\lambda)}$ is
minimized (among partitions $\mu,\lambda$ of size $n$) by
$\mu=(n),\lambda=(1^n)$ or $\mu=(1^n),\lambda=(n)$.
\end{prop}

\begin{proof} Proposition \ref{pos} gives that \[ \frac{(U
D)^r [\mu,\lambda]}{\pi(\lambda)} = \sum_{k=0}^n A_n(r,k) \frac{ (U^k
D^k)[\mu,\lambda] }{\pi(\lambda)},\] with all $A_n(r,k) \geq 0$. The
proposition now follows from three observations:
\begin{itemize}
\item All terms in the sum are non-negative. Indeed, Proposition \ref{pos} gives
that all $A_n(r,k) \geq 0$, and $U,D$ were defined probabilistically.
\item If $\mu=(n),\lambda=(1^n)$ or $\mu=(1^n),\lambda=(n)$, then the summands for
$0 \leq k \leq n-2$ vanish. Indeed, for such $k$ it is impossible to move
from the partition $\mu$ to the partition $\lambda$ by removing $k$ boxes
one at a time and then reattaching $k$ boxes one at a time.
\item The $k=n-1$ and
$k=n$ summands are each independent of both $\mu$ and $\lambda$. Indeed,
$D^{n-1}(\mu)$ is equal to $(1)$ for any partition $\mu$ of size $n$.
Since the up chain preserves the Pitman distribution, it follows that
$U^{n-1} [(1),\lambda]=\pi(\lambda)$, so that the $k=n-1$ summand is
$A_n(r,n-1)$. Similarly, the $k=n$ summand is $A_n(r,n)$.
\end{itemize} \end{proof}

The following corollary will be helpful.

\begin{cor} \label{sepdist} Consider down-up walk with parameters $\theta>0$
and $0 \leq \alpha <1$ on the nth level of the Kingman lattice. Then
$s^*(r)=\pp(T>r)$ where $T$ is the sum of independent geometrics with
parameters $\frac{i(\theta+i-1)}{n(\theta+n-1)}$ for $2 \leq i \leq n$.
\end{cor}

\begin{proof} By Proposition \ref{sepking}, $s^*(r) = 1 - \frac{(UD)^r((n),(1^n))}
{\pi(1^n)}$. By Corollary \ref{eigenking}, the down-up walk has $n$
distinct eigenvalues, namely $1$ and $1-\frac{i
(\theta+i-1)}{n(\theta+n-1)}$ for $2 \leq i \leq n$. Since the distance
between $(n)$ and $(1^n)$ is $n-1$, the result follows from Propositions
\ref{noeigenvec} and \ref{geomeig}.
\end{proof}

Theorem \ref{exactking} gives the precise asymptotic behavior of $s^*(r)$
in the special case that $\theta=1$.

\begin{theorem} \label{exactking} Let $s^*(r)$ be the maximal separation
distance after $r$ iterations of down-up walk on the Kingman lattice, in
the special case that $\theta=1$ and $0 \leq \alpha <1$.
\begin{enumerate}
\item \[ s^*(r) = 2 \sum_{i=2}^n (-1)^i (i^2-1) \frac{(n!)^2}{(n-i)!(n+i)!}
\left( 1-\frac{i^2}{n^2} \right)^r.\]
\item For $c>0$ fixed,
\[ \lim_{n \rightarrow \infty} s^*(cn^2) = 2 \sum_{i=2}^{\infty}
 (-1)^i (i^2-1) e^{-ci^2}.\]
\end{enumerate}
\end{theorem}

\begin{proof} By Proposition \ref{sepking}, one has that \[ s^*(r) = 1 -
\frac{(UD)^r[(n),(1^n)]}{\pi(1^n)}.\] By Corollary \ref{eigenking}, the
chain has $n$ distinct eigenvalues. Since the distance between $(n)$ and
$(1^n)$ is $n-1$, it follows from Proposition \ref{noeigenvec} that
\begin{eqnarray*} s^*(r) & = & \sum_{i=2}^n \left( 1 - \frac{i^2}{n^2}
\right)^r \prod_{2 \leq j \leq n \atop j \neq i} \frac{
\frac{j^2}{n^2}}{\frac{j^2}{n^2}-\frac{i^2}{n^2}}\\
& = & \sum_{i=2}^n \left( 1 - \frac{i^2}{n^2} \right)^r \prod_{2 \leq j
\leq n \atop j \neq i} \frac{j^2}{(j-i)(j+i)}, \end{eqnarray*} and the
first assertion follows by elementary simplifications.

For part 2 of the theorem, we claim that for $c>0$ fixed there is a
constant $i_c$ (depending on $c$ but not $n$) such that for $i \geq i_c$,
the summands in \[ 2 \sum_{i=2}^n (-1)^i (i^2-1)
\frac{(n!)^2}{(n-i)!(n+i)!} \left( 1-\frac{i^2}{n^2} \right)^{cn^2}\] are
decreasing in magnitude (and alternating in sign). Part 2 of the theorem
follows from this claim, since then one can take limits for each fixed
$i$. To prove the claim, note that the summands are decreasing in
magnitude if $i \geq \sqrt{n}$, since one checks that $(i^2-1)
\frac{(n!)^2}{(n-i)!(n+i)!}$ is a decreasing function of $i$ when $i \geq
\sqrt{n}$. Since $\frac{(n!)^2}{(n-i)!(n+i)!}$ is a decreasing function of
$i$, to handle $i \leq \sqrt{n}$ one needs only to show that
\[ \frac{i^2-1}{(i+1)^2-1} e^{cn^2 [\log(1-i^2/n^2) - \log(1-(i+1)^2/n^2)] }
> 1 \] for $i \geq i_c$, a constant depending on $c$ but not $n$. Using
that $\log(1-x) \geq -x -x^2$ for $0<x<\frac{1}{2}$ and that $\log(1-x)
\leq -x$ for $0 < x <1$, one has that \[ cn^2 \left[ \log(1-i^2/n^2) -
\log(1-(i+1)^2/n^2) \right] \geq c (i+1)^2 - ci^2 - c \frac{i^4}{n^2} \geq
2ic
\] since $i \leq \sqrt{n}$. Clearly $\frac{i^2-1}{(i+1)^2-1} e^{2ic} > 1 $
for $i$ large enough, completing the proof. \end{proof}

For general values of $\theta$, we do not have a result as precise as
Theorem \ref{exactking}, but obtain explicit upper and lower bounds for
the separation distance mixing time. Note that when $\theta$ is fixed and
$n$ is growing, these bounds are of order $n^2$.

\begin{cor} Let $n^*_{\frac{1}{2}}$ be the separation distance mixing
time for down-up walk (with parameters $\theta>0$ and $0 \leq \alpha<1$)
on the nth level of Kingman's lattice. Then $\frac{\ee[T]}{2} \leq
n^*_{\frac{1}{2}} \leq 2 \ee[T]$, where $T$ is as in Corollary
\ref{sepdist}. Moreover if $\theta=1$ then \[ \ee[T] = \sum_{i=2}^n
\frac{n^2}{i^2} \sim n^2 \left( \frac{\pi^2}{6}-1 \right),\] and if
$\theta \neq 1$ then \begin{eqnarray*} \frac{n(\theta+n-1)}{\theta-1} \log
\left( \frac{(n+1)(\theta+1)}{2(n+\theta)} \right) \leq \ee[T] =
\sum_{i=2}^n \frac{n(\theta+n-1)}{i(\theta+i-1)} \\ \leq
\frac{n(\theta+n-1)}{\theta-1} \log \left( \frac{n \theta}{n+\theta-1}
\right). \end{eqnarray*}
\end{cor}

\begin{proof} Lemma \ref{tbound} gives that $\frac{\ee[T]}{2} \leq
n^*_{\frac{1}{2}} \leq 2 \ee[T]$ and Corollary \ref{sepdist} gives that
$\ee[T] = \sum_{i=2}^n \frac{n(\theta+n-1)}{i(\theta+i-1)}$. To complete
the proof of the upper bound, note that \[ \sum_{i=2}^n
\frac{1}{i(\theta+i-1)} \leq \int_1^n \frac{1}{t(\theta+t-1)} dt =
 \frac{1}{\theta-1} \log \left( \frac{n \theta}{n+\theta-1}
\right). \] For the lower bound, note that \[ \sum_{i=2}^n
\frac{1}{i(\theta+i-1)} \geq \int_2^{n+1} \frac{1}{t(\theta+t-1)} dt =
\frac{1}{\theta-1} \log \left( \frac{(n+1)(\theta+1)}{2(n+\theta)}
\right). \] \end{proof}

To conclude, we relate separation distance of the up-down chain to
separation distance of the down-up chain.

\begin{prop} \label{samerateking} Let $s^*_{UD_n}(r)$ be the maximal separation
distance after $r$ iterations of the down-up chain (with parameters
$\theta>0$ and $0 \leq \alpha<1$) on the Kingman lattice, and let
$s^*_{DU_n}(r)$ be the corresponding quantity for the up-down chain. Then
\[ s^*_{DU_n}(r) = s^*_{UD_{n+1}}(r+1) \] for all $n \geq 1, r \geq 0$.
\end{prop}

\begin{proof} The method is the same as for Proposition \ref{samerate}.
The eigenvalues of $DU_n$ are $1$ and $t_j:=1-\frac{j(\theta+j-1)}
{(n+1)(\theta+n)}$ (for $2 \leq j \leq n$) yielding that \[ s^*_{DU_n}(r)
= \sum_{j=2}^n (t_j)^r  \prod_{k \neq j \atop 2 \leq k \leq n} \left(
\frac{1-t_k}{t_j-t_k} \right).\] The eigenvalues of $UD_{n+1}$ are $1$ and
$t_j$ (for $2 \leq j \leq n+1$) yielding that \[ s^*_{UD_{n+1}}(r+1) =
\sum_{j=2}^{n+1} (t_j)^{r+1} \prod_{k \neq j \atop 2 \leq k \leq n+1}
\left( \frac{1-t_k}{t_j-t_k} \right).\] The result follows since
$t_{n+1}=0$.
\end{proof}

\section{Other examples} \label{other}

This section treats other examples to which the commutation relation
methodology applies. After discussing two classical examples
(Bernoulli-Laplace models and subspace walks), we determine precise
separation distance asymptotics for a non-standard hypercube example.

We focus on the down-up chain but for readers interested in the up-down
chain mention the relation $s^*_{DU_n}(r)=s^*_{UD_{n+1}}(r+1)$ (which is
true for the same reasons as in the Young and Kingman examples). This
holds for all examples in this section except for the subset walk on
$\lfloor \frac{n}{2} \rfloor$ sets or the subspace walk on $\lfloor
\frac{n}{2} \rfloor$ spaces (in these exceptional cases the two chains
have the same separation distance asymptotics).

\subsection{Bernoulli-Laplace models} We analyze random walk on size
$j$ subsets of an $n$ element set, where $0< 2j \leq n$. From a subset $S$
of size $j$, a step proceeds by first removing one of the $j$ elements
uniformly at random, and then randomly adding in one of the $n-j+1$
elements in $S-j$. The stationary distribution is the uniform distribution
on subsets of size $j$. This chain appears when analyzing the
Bernoulli-Laplace model, in which there are two urns, the left containing
$j$ red balls, the right containing $n-j$ black balls, and at each step a
ball is picked uniformly at random in each urn, and the two balls are
switched.

It will be useful to let $P$ be the Boolean lattice of rank $n$; the
elements of $P$ are the subsets of $\{1,\cdots,n\}$ and $S \leq T$ in the
partial order if $S \subseteq T$. Letting $U,D$ be the up and down
operators for this poset, Stanley \cite{St3} observed that \[ D_{j+1} U_j
= U_{j-1} D_j + (n-2j) I_j,\] for  $0 \leq j \leq n$. For our purposes, it
is more convenient to work with the normalized operators \[ \tilde{U}_j =
\frac{1}{n-j} U_j \ , \ \tilde{D}_j=\frac{1}{j} D_j.\] Then the random
walk on size $j$ subsets of $\{1,\cdots,n\}$ is given by the operator
$\tilde{U} \tilde{D}_j$. Stanley's commutation relation becomes \[
\tilde{D}_{j+1} \tilde{U}_j = a_j \tilde{U}_{j-1} \tilde{D}_j + (1-a_j)
I_j \] with $a_j=\frac{j(n-j+1)}{(j+1)(n-j)}$.

As a consequence of Theorem \ref{eigenvalgen}, one obtains the eigenvalues
of $\tilde{U} \tilde{D}_j$. This goes back at least to Karlin and McGregor
\cite{KM}.

\begin{cor} \label{berneig} The eigenvalues of $\tilde{U}\tilde{D}_{j}$
are $$\left\{ \begin{array}{ll} 1 & \mbox{multiplicity} \ 1 \\
1-\frac{i(n-i+1)}{j(n-j+1)} & \mbox{multiplicity}
 \ {n \choose i} - {n \choose i-1} \ (1 \leq i \leq j) \end{array} \right.$$
\end{cor}

Proposition \ref{sepatt} gives information about separation distance. The
proof in \cite{DF} used the theory of birth-death chains, and the fact
that the Bernoulli-Laplace chain can be reduced to a birth death chain
(look at the number of red balls in the right urn). Our proof uses
commutation relations.

\begin{prop} (\cite{DF}) \label{sepatt} Consider the random walk $\tilde{U} \tilde{D}_j$
on size $j$ subsets of $\{1,\cdots,n\}$. Let $r$ be a non-negative
integer, and let $\pi$ be the uniform distribution on $j$ element subsets
of $\{1,\cdots,n\}$.
\begin{enumerate}
\item The quantity $\frac{(\tilde{U} \tilde{D})^r[S,T]}{\pi(T)}$ is
minimized (among pairs of $j$ element subsets of $\{1,\cdots,n\}$) by any
$S,T$ such that $S \cap T = \emptyset$.
 \item \[ s^*(r) = \pp(X>r), \]
where $X$ is the sum of independent geometrics having parameters
$\frac{i(n-i+1)}{j(n-j+1)}$ for $1 \leq i \leq j$.
\end{enumerate}
\end{prop}

\begin{proof} Given a linear operator $A:\mathbb{C}P \mapsto \mathbb{C}P$,
and subsets $S,T$ of $\{1,\cdots,n\}$ of size $j$, let $A[S,T]$ denote the
coefficient of $T$ in $A(S)$. Proposition \ref{pos} gives that \[
\frac{(\tilde{U} \tilde{D})^r[S,T]}{\pi(T)} = \sum_{k=0}^j A_j(r,k)
\frac{\tilde{U}^k \tilde{D}^k[S,T]}{\pi(T)}, \] with all $A_j(r,k) \geq
0$. The first part of the proposition now follows from three observations:
\begin{itemize}
\item All terms in the sum are non-negative. Indeed, all $A_j(r,k) \geq 0$
and $\tilde{U},\tilde{D}$ were defined probabilistically.
\item If $S \cap T = \emptyset$, then the summands for $0 \leq k \leq j-1$
 all vanish. This is clear since for such $k$, $\tilde{U}^k
 \tilde{D}^k[S,T]=0$.

\item The $k=j$ summand is independent of both $S$ and $T$. Indeed,
$\tilde{D}^{j}(S)=\emptyset$ for any $S$ of size $j$, and
$\tilde{U}^j(\emptyset)$ is uniformly distributed among the size j subsets
of $\{1,\cdots,n\}$. Hence the $k=j$ summand is equal to $A_j(r,j)$.
\end{itemize}

For the second part of the proposition, Corollary \ref{berneig} gives that
$\tilde{U} \tilde{D}_j$ has $j+1$ distinct eigenvalues. Letting $x=S,y=T$
where $S \cap T = \emptyset$, one has that $dist(x,y)=j$. The result now
follows from Propositions \ref{noeigenvec} and \ref{geomeig}. \end{proof}

In fact there is another proof of part 2 of Proposition \ref{sepatt} which
uses only combinatorial properties of the sequence $A_j(r,j)$.

\begin{proof} (Second proof of part 2 of Proposition \ref{sepatt}) The
proof of part 1 of Proposition \ref{sepatt} gives that
$s^*(r)=1-A_j(r,j)$, where $A_j(r,j)$ is defined in Proposition \ref{pos}.
Letting $[x^n] f(x)$ denote the coefficient of $x^n$ in a power series
$f(x)$, Proposition \ref{gen} gives that
\begin{eqnarray*}
A_j(r,j) & = & [x^r] \frac{x^j \prod_{i=1}^j
\frac{(j-i+1)(n-j+i)}{j(n-j+1)}}{\prod_{i=1}^j 1 - x \left( 1-
\frac{(j-i)(n-j+i+1)}{j(n-j+1)} \right)} \\
& = & [x^r] \frac{1}{1-x} \prod_{i=1}^j  \frac{x
\frac{(j-i+1)(n-j+i)}{j(n-j+1)}}{1-x \left( 1 -
\frac{(j-i+1)(n-j+i)}{j(n-j+1)} \right)} \\
& = & [x^r] \frac{1}{1-x} \prod_{i=1}^j  \frac{x
\frac{i(n-i+1)}{j(n-j+1)}}{1-x \left( 1 - \frac{i(n-i+1)}{j(n-j+1)}
\right)}. \end{eqnarray*} The last step used the change of variables $i
\mapsto j+1-i$.

Note that if $Z$ is geometric with parameter $p$, then $Z$ has probability
generating function
\[ \sum_{i \geq 0} x^i \pp(Z=i) = \frac{xp}{1-x(1-p)} .\] Thus $A_j(r,j)$
is the probability that the convolution of geometrics with parameters
$\frac{i(n-i+1)}{j(n-j+1)}$ is at most $r$, and the result follows.
\end{proof}

The asymptotic behavior of $s^*(r)$ (in continuous time) is studied in
detail in \cite{DSa}, using a continuous time analog of part 2 of
Proposition \ref{sepatt} (in which geometrics are replaced by
exponentials). A similar analysis can be carried out in discrete time. For
instance if $j \leq \frac{n}{2}$ tends to infinity, there is a separation
cutoff at time $t_{n,j}=\frac{j(n-j)}{n} \log(j)$. For information
concerning convergence in the total variation metric, see \cite{Be} or
\cite{DSh}.

\subsection{Subspace walks}

This is a $q$-analog of the previous example. The random walk is on
j-dimensional subspaces of an n-dimensional vector space over a finite
field $\mathbb{F}_q$, where $0<2j \leq n$. From a j-dimensional subspace
$S$, a step of the walk proceeds by first choosing uniformly at random a
$j-1$ dimensional subspace $W$ contained in $S$, and then choosing
uniformly at random a $j$ dimensional subspace $T$ containing $W$.

Up to holding, this random walk is equivalent to the nearest neighbor walk
on the graph of $j$ dimensional subspaces, where two subspaces are
connected by an edge if their intersection has dimension $j-1$. As
discussed in \cite{Be}, \cite{D'A}, the eigenvalues of this walk are known
and sharp total variation distance estimates can be obtained by studying a
related birth-death chain on $\{0,\cdots,j\}$, which is just the
associated graph distance process.

To revisit this example using commutation relations, let $P$ be the
subspace lattice of an $n$-dimensional vector space over a finite field
$\mathbb{F}_q$. Letting $U,D$ be the up and down operators for the poset
$P$, Stanley \cite{St3} observed that \[ D_{j+1} U_j = U_{j-1} D_j +
\left( \frac{q^{n-j}-1}{q-1} - \frac{q^{j}-1}{q-1} \right) I_j, \] for $0
\leq j \leq n$. For our purposes it is convenient to renormalize the
operators as
\[ \tilde{U}_j = \frac{q-1}{q^{n-j}-1} U_j \ , \ \tilde{D}_j =
\frac{q-1}{q^j-1} D_j.\] Then the random walk on $j$ dimensional subspaces
is given by $\tilde{U} \tilde{D}_j$, and one checks that the commutation
relation becomes \[ \tilde{D}_{j+1} \tilde{U}_j = a_j \tilde{U}_{j-1}
\tilde{D}_j + (1-a_j) I_j
\] where $a_j=\frac{(q^{n-j+1}-1)(q^j-1)}{(q^{n-j}-1)(q^{j+1}-1)}$.

As an immediate consequence of this commutation relation and Theorem
\ref{eigenvalgen}, one obtains the eigenvalues of the subspace walk.

\begin{cor} \label{berneigq} The eigenvalues of $\tilde{U} \tilde{D}_j$
 are $$\left\{ \begin{array}{ll} 1 & \mbox{multiplicity} \ 1 \\
1-\frac{(q^{n-i+1}-1)(q^i-1)}{(q^{n-j+1}-1)(q^j-1)} &
\mbox{multiplicity} \ \qb{n}{i} - \qb{n}{i-1} \ (1 \leq i \leq j)
\end{array} \right.$$ Here $\qb{n}{i}$ denotes the number of $i$-dimensional
subspaces of an n-dimensional vector space over $\mathbb{F}_q$.
\end{cor}

Proposition \ref{sepq} gives a result about separation distance. This also
follows from the birth-death chain theory in \cite{DF}.

\begin{prop} \label{sepq} Consider the random walk $\tilde{U}\tilde{D}_j$
on j-dimensional subspaces of an $n$ dimensional vector space $V$ over
$\mathbb{F}_q$. Let $r$ be a non-negative integer, and let $\pi$ be the
uniform distribution on $j$-dimensional subspaces of $V$.
\begin{enumerate}
\item The quantity $\frac{(\tilde{U} \tilde{D})^r[S,T]}{\pi(T)}$ is
minimized (among pairs of $j$ dimensional subspaces of $V$) by any $S,T$
such that $S \cap T = 0$.
\item One has that $s^*(r) =
\pp(X>r)$, where $X$ is the sum of independent geometrics with parameters
$\frac{(q^{n-i+1}-1)(q^i-1)}{(q^{n-j+1}-1)(q^{j}-1)}$, for $1 \leq i \leq
j$.
\end{enumerate}

\end{prop}

\begin{proof} The proof method for both parts is the same as for the proof
of Proposition \ref{sepatt}; one need only replace the word ``subset'' by
``subspace'' and the word ``size'' by ``dimension''. Note that the second
proof of part of Proposition \ref{sepatt} also carries over to the
subspace setting.
\end{proof}

Concerning the asymptotic behavior of $s^*(r)$, we note that \cite{DSa}
gives results (in the continuous time case), using an analog of part 2 of
Proposition \ref{sepq} in which the geometrics are replaced by
exponentials. Their method can be transferred to the discrete time
setting. For instance if $j \leq \frac{n}{2}$ tends to infinity, there is
a separation cutoff at time $t_{n,j}=j$.

\subsection{Gibbs sampler for hypercube}

The main object of study in this example is the birth-death chain on the
set $\{0,1,\cdots,n\}$ with transition probabilities \[ K(x,x-1) =
\frac{x}{n} (1-p) , \ K(x,x) = \frac{x}{n} p + \left( 1-\frac{x}{n}
\right) (1-p)\] \[  K(x,x+1) = p \left( 1-\frac{x}{n} \right).\] Here $0 <
p < 1$ and the stationary distribution of this chain is the p-binomial
distribution $\pi(x) = {n \choose x} p^x (1-p)^{n-x}$.

We remark that this Markov chain is the distance chain for the Gibbs
sampler on the hypercube, used to sample from the distribution in which a
length $n$ 0-1 vector is assigned probability $p^x (1-p)^{n-x}$, where $x$
is the number of $1$'s in the vector. For general $p$ we have not seen
this exact analyzed chain in the literature (though possibly it has been
studied). Different birth-death chains with the same stationary
distribution are studied as examples in \cite{DSa}. Our birth-death chain
has the property that the eigenvalues are independent of $p$ (see
Corollary \ref{eigenpas}); the examples in \cite{DSa} do not.

To motivate the definition of up and down operators, we note that the
birth-death chain in this section is, in the terminology of Section
\ref{downup}, an example of a down-up Markov chain. The poset we use is
Pascal's lattice: the vertices of the nth level are labeled by pairs
$(x,n)$ where $x=0,1,\cdots,n$. The only edges are $(x,n) \nearrow
(x,n+1)$ and $(x,n) \nearrow (x+1,n+1)$, each with multiplicity 1. Then
the dimension of the vertex $(x,n)$ is ${n \choose x}$. One checks that
the probability distributions $M_n((x,n))= {n \choose x} p^x (1-p)^{n-x}$
are coherent with respect to Pascal's lattice \cite{K}, and computes that
the corresponding up and down chains are given by
\[ U_n [(x,n)] = (1-p)\cdot (x,n+1) + p \cdot (x+1,n+1) \]
\[ D_n [(x,n)] = \left( 1 - \frac{x}{n} \right) \cdot (x,n-1) +
\frac{x}{n} \cdot (x-1,n-1).\] From this one sees that our birth-death
chain is precisely the down-up chain $UD_n$ on Pascal's lattice.

\begin{prop} \label{relevpas} Letting $a_n = \frac{n}{n+1}$, one has that
\[ D_{n+1} U_n = a_n U_{n-1} D_n + (1-a_n) I_n.\]
\end{prop}

\begin{proof} This is straightforward to check from the definitions of
$U$ and $D$. \end{proof}

Corollary \ref{eigenpas} determines the eigenvalues of the down-up walk on
Pascal's lattice. It is curious that they are independent of $p$.

\begin{cor} \label{eigenpas} The eigenvalues of $U D_n$ are $1-\frac{i}{n}$
with multiplicity 1, for $0 \leq i \leq n$.
\end{cor}

\begin{proof} This is immediate from Theorem \ref{eigenvalgen} and
Proposition \ref{relevpas}. \end{proof}

Proposition \ref{seppas} determines where the maximal separation distance
is attained.

\begin{prop} \label{seppas} Let $\pi$ be the p-binomial distribution and
let $r$ be a non-negative integer. The quantity $\frac{ (U
D)^r[(x,n),(y,n)]}{\pi((y,n))}$ is minimized (among $0 \leq x,y \leq n$)
by $x=0,y=n$ or $x=n,y=0$.
\end{prop}

\begin{proof} Given a linear operator $B: \mathbb{C}P_n \mapsto \mathbb{C}P_n$,
let $B[(x,n),(y,n)]$ denote the coefficient of $(y,n)$ in $B(x,n)$.
Proposition \ref{pos} gives that \[ \frac{(U D)^r[(x,n),(y,n)]}{\pi(y,n)}
= \sum_{k=0}^n A_n(r,k) \frac{U^k D^k [(x,n),(y,n)]}{\pi(y,n)}, \] with
all $A_n(r,k) \geq 0$. The proposition now follows from three facts:
\begin{itemize}
\item All terms in the sum are non-negative. Indeed, all $A_n(r,k) \geq 0$
and $U,D$ were defined probabilistically.
\item If $x=0,y=n$ or $x=n,y=0$, the summands for $0 \leq k \leq n-1$ all vanish.
\item The $k=n$ summand is independent of both $x$ and $y$. Indeed,
$D^n(x,n)=(0,0)$ and the coefficient of $(y,n)$ in $U^n(0,0)$ is
$\pi(y,n)$. So the $k=n$ summand is exactly $A_n(r,n)$. \end{itemize}
\end{proof}

Finally, we determine the exact asymptotic behavior of $s^*(r)$ for this
example.

\begin{prop} \label{cutpas} Consider the random walk $UD_n$ corresponding to
the $p$-binomial distribution. Let $r$ be a non-negative integer.
\begin{enumerate}
\item $s^*(r) = \pp(X>r)$ where $X$ is the sum of independent geometrics with
parameters $\frac{i}{n}$ for $1 \leq i \leq n$.
\item $s^*(r) = 1 - \frac{n!
S(r,n)}{n^r}$ where $S(r,k)$ is a Stirling number of the second kind (i.e.
the number of partitions of an r set into k blocks).
\item For $c$ fixed in
$\mathbb{R}$ and $n \rightarrow \infty$, \[ s^*(n \log(n)+cn) = 1 -
e^{-e^{-c}} + O \left( \frac{\log(n)}{n} \right) .\]
\end{enumerate}
\end{prop}

\begin{proof} Proposition \ref{seppas} gives that $s^*(r)=1 - \frac{(UD)^r((0,n),(y,n))}
{\pi(y,n)}$. By Corollary \ref{eigenpas} the chain has $n+1$ distinct
eigenvalues. Hence the first assertion follows from Proposition
\ref{noeigenvec} (with $x=(0,n)$ and $y=(n,n)$), and Proposition
\ref{geomeig}.

For the second assertion, it follows from the proof of Proposition
\ref{seppas} and Proposition \ref{pos} that $s^*(r)=1-A_n(r,n)$ where
$A_n(r,k)$ satisfies the recurrence \[ A_n(r,k) = \frac{n-k+1}{n}
A_n(r-1,k-1) + \frac{k}{n} A_n(r-1,k) \] with initial condition
$A_n(0,m)=\delta_{0,m}$. It is straightforward to check that
$A_n(r,k)=\frac{n! S(r,k)}{n^r (n-k)!}$ solves the recurrence, using the
recurrence for Stirling numbers \[ S(r,k) = S(r-1,k-1) + kS(r-1,k)  \] on
page 33 of \cite{St1}.

For the third assertion, it follows from the second assertion and the
argument in part 2 of Theorem \ref{exactplan} that $s^*(r)=1-P(n,r,n)$,
where $P(n,r,n)$ is the probability of $n$ occupied boxes when $r$ balls
are dropped into $n$ boxes. The result now follows from asymptotics of the
coupon collector's problem, as in the proof of Theorem \ref{exactplan}.
\end{proof}

{\it Remark:} The waiting time for $n$ boxes to all be occupied when balls
are randomly dropped into them one at a time is a convolution of
independent geometrics with parameters $\frac{i}{n}$ for $1 \leq i \leq
n$. Thus part 3 of Proposition \ref{cutpas} can be proved without using
part 2 of Proposition \ref{cutpas}. Our reason for using part 2 was to
illustrate that one can sometimes usefully solve the recursion for the
combinatorially defined quantities $A_n(r,k)$.

\section*{Acknowledgements} The author received funding from NSF grant
DMS-0503901, and thanks the referee for helpful comments.


\begin{thebibliography}{AAA}

\bibitem [AD1]{AD1} Aldous, D. and Diaconis, P., Shuffling cards and
stopping times, {\it Amer. Math. Monthly} {\bf 93} (1986), 333-348.

\bibitem [AD2]{AD2} Aldous, D. and Diaconis, P., Strong uniform times and
finite random walks, {\it Adv. in Appl. Math.} {\bf 8} (1987), 69-97.

\bibitem [An]{An} Andrews, G., {\it The theory of partitions}, Cambridge
University Press, Cambridge, 1984.

\bibitem [Be]{Be} Belsley, E., Rates of convergence of random walk on
distance regular graphs, {\it Probab. Theory Relat. Fields} {\bf 112}
(1998), 493-533.

\bibitem [B]{B} Borodin, A., Multiplicative central measures in the Schur
graph, in {\it Representation theory, dynamical systems, combinatorial and
algorithmic methods II} (A.M. Vershik, ed.), Zap. Nauchn. Sem. POMI {\bf
240}, Nauka, St. Petersburg, 1997, 44-52 (Russian); English translation in
{\it J. Math. Sci. (New York)} {\bf 96} (1999), 3472-3477.

\bibitem [BOO] {BOO} Borodin, A., Okounkov, A., and Olshanski, G.,
Asymptotics of Plancherel measures for symmetric groups, {\it J. Amer.
Math. Soc.} {\bf 13} (2000), 481-515.

\bibitem [BO1] {BO1} Borodin, A. and Olshanski, G., Harmonic functions
on multiplicative graphs and interpolation polynomials, {\it
Electron. J. Combin.} {\bf 7} (2000), Research paper 28, 39 pages
(electronic).

\bibitem [BO2]{BO2} Borodin, A. and Olshanski, G., Infinite
dimensional diffusions as limits of random walks on partitions, arXiv:
math.PR/0706.1034 (2007).

\bibitem [BO3]{BO3} Borodin, A. and Olshanski, G., Markov processes on
partitions, {\it Probab. Theory Relat. Fields} {\bf 135} (2006), 84-152.

\bibitem [BO4]{BO4} Borodin, A. and Olshanski, G., Z-measures on
partitions and their scaling limits, {\it European J. Combin.} {\bf 26}
(2005), 795-834.

\bibitem [Br]{Br} Brown, M., Spectral analysis, without eigenvectors,
for Markov chains, {\it Probab. Eng. Inform. Sci.} {\bf 5} (1991),
131-144.

\bibitem [BS]{BS} Brown, M. and Shao, Y., Identifying coefficients in the
spectral representation for first passage time distributions, {\it Probab.
Eng. Inform. Sci.} {\bf 1} (1987), 69-74.

\bibitem [C]{C} Chatterjee, S., Stein's method for concentration
inequalities, {\it Probab. Theory Relat. Fields} {\bf 138} (2007),
305-312.

\bibitem [CDM]{CDM} Chatterjee, S., Diaconis, P., and Meckes, E.,
Exchangeable pairs and Poisson approximation, {\it Probab. Surv.} {\bf 2}
(2005), 64-106.

\bibitem [D'A]{D'A} D'Aristotle, A., The nearest neighbor random walk
on subspaces of a vector space and rate of convergence, {\it J. Theoret.
Probab.} {\bf 8} (1993), 321-346.

\bibitem [D]{D} Diaconis, P.,  The cutoff phenomenon in finite Markov
chains, {\it Proc. Nat. Acad. Sci. U.S.A.} {\bf 93} (1996), 1659-1664.

\bibitem [DF]{DF} Diaconis, P. and Fill, J., Strong stationary times via a
new form of duality, {\it Ann. Probab.} {\bf 18} (1990), 1483-1522.

\bibitem [DH]{DH} Diaconis, P. and Hanlon, P., Eigen-analysis for some
examples of the Metropolis algorithm, in {\it Hypergeometric functions on
domains of positivity, Jack polynomials, and applications}, 99-117,
Contemp. Math. 138, 1992.

\bibitem [DSa]{DSa} Diaconis, P. and Saloff-Coste, L., Separation cutoffs
for birth death chains, {\it Ann. Appl. Probab.} {\bf 16} (2006),
2098-2122.

\bibitem [DSh]{DSh} Diaconis, P. and Shahshahani, M., Time to reach
stationarity in the Bernoulli-Laplace diffusion model, {\it SIAM
J. Math. Anal.} {\bf 18} (1987), 208-218.

\bibitem [Ew]{Ew} Ewens, W.J., Population genetics theory: the past and the
future, in {\it Mathematical and statistical developments of evolutionary
theory}. Kluwer, Dordrecht, 1990, 117-228.

\bibitem [Fo]{Fo} Fomin, S., Duality of graded graphs, {\it J. Algebraic
Combin.} {\bf 3} (1994), 357-404.

\bibitem [F1]{F1} Fulman, J., Stein's method and Plancherel measure of
the symmetric group, {\it Trans. Amer. Math. Soc.} {\bf 357} (2005),
555-570.

\bibitem [F2]{F2} Fulman, J., Convergence rates of random walk on
irreducible representations of finite groups, {\it
J. Theoret. Probab.}, to appear.

\bibitem [F3]{F3} Fulman, J., Separation cutoffs for random walk on
irreducible representations, arXiv: math.PR/0703291 (2007).

\bibitem [F4]{F4} Fulman, J., Stein's method and random character ratios,
{\it Transac. Amer. Math. Soc.}, to appear.

\bibitem [F5]{F5} Fulman, J., Stein's method, Jack measure, and the
Metropolis algorithm, {\it J. Combin. Theory Ser. A.} {\bf 108} (2004),
275-296.

\bibitem [HH]{HH} Hoffman, P. and Humphreys, J., {\it Projective
representations of the symmetric group}, Oxford University Press, New
York, 1992.

\bibitem [I]{I} Ivanov, V., Plancherel measure on shifted Young diagrams,
in {\it Representation theory, dynamical systems, and asymptotic
combinatorics}, Amer. Math. Soc. Transl. Ser. 2, 217, (2006), 73-86.

\bibitem [KM]{KM} Karlin, S. and McGregor, J., Ehrenfest urn models,
{\it J. Appl. Probab.} {\bf 2} (1965), 352-376.

\bibitem [K]{K} Kerov, S., The boundary of Young lattice and random Young
tableaux, {\it Formal power series and algebraic combinatorics}, DIMACS
Ser. Discrete Math. Theoret. Comput. Sci. {\bf 24}, Amer. Math. Soc.,
Providence, RI, (1996), 133-158.

\bibitem [KOV1]{KOV1} Kerov, S., Olshanski, G., and Vershik, A.,
Harmonic analysis on the infinite symmetric group. A deformation of the
regular representation, {\it C.R. Acad. Sci. Paris S\'{e}r. I Math.} {\bf
316} (1993), 773-778.

\bibitem [KOV2]{KOV2} Kerov, S., Olshanski, G., and Vershik, A., Harmonic
analysis on the infinite symmetric group, {\it Invent. Math.} {\bf 158}
(2004), 551-642.

\bibitem [Mac]{Mac} Macdonald, I., {\it Symmetric functions and Hall
polynomials}, Second edition, Oxford University Press, New York, 1995.

\bibitem [Mat]{Mat} Matsumoto, S., Correlation functions of the shifted
Schur measure, {\it J. Math. Soc. Japan} {\bf 57} (2005), 619-637.

\bibitem [O1]{O1} Okounkov, A., $SL(2)$ and $z$ measures, in {\it Random
matrix models and their applications}, 407-420, Math. Sci. Res. Inst.
Publ. 40, Cambridge Univ. Press, Cambridge, 2001.

\bibitem [O2]{O2} Okounkov, A., The uses of random partitions, in {\it
XIVth International Congress on Mathematical Physics}, 379-403, World Sci.
Publ., Hackensack, NJ, 2005.

\bibitem [Pa]{Pa} Pak, I., {Random walk on groups: strong uniform time
approach}, Ph.D. Thesis, Harvard University, 1997.

\bibitem [Pe]{Pe} Petrov, L., Two-parameter family of diffusion processes in
the Kingman simplex, arXiv: math.PR/0708.1930 (2007).

\bibitem [R]{R} Ross, N., Step size in Stein's method of exchangeable pairs,
preprint, 2007.

\bibitem[Sag]{Sag} Sagan, B., {\it The symmetric group. Representations,
combinatorial algorithms, and symmetric functions}, Springer-Verlag, New
York, 1991.

\bibitem[Sal]{Sal} Saloff-Coste, L., Random walk on finite groups, in
{\it Probability on discrete structures}, 263-346, Encyclopedia Math. Sci.
110, Springer, Berlin, 2004.

\bibitem[St1]{St1} Stanley, R., {\it Enumerative combinatorics,
Vol. 1}, Wadsworth \& Brooks/Cole, Monterey, 1986.

\bibitem[St2]{St2} Stanley, R., Differential posets, {\it J. Amer. Math.
Soc.} {\bf 1} (1988), 919-961.

\bibitem[St3]{St3} Stanley, R., Variations on differential posets, in {\it
Invariant theory and tableaux}, IMA Vol. Math. Appl. 19, Springer, New
York, 1990, 145-165.

\bibitem [TW]{TW} Tracy, C. and Widom, H., A limit theorem for shifted
Schur measures, {\it Duke Math. J.} {\bf 123} (2004), 171-208.

\bibitem [Wi]{Wi} Wilf, H., {\it Generatingfunctionology}, Second edition.
Academic Press, Inc., Boston, 1994.

\bibitem [Wk]{Wk} Wilkinson, J., {\it The algebraic eigenvalue problem}, Oxford
University Press, Oxford, 1988.
\end{thebibliography}
\end{document}